\documentclass[preprint,11pt]{amsart}

\usepackage{ulem}
\usepackage{graphicx}
\usepackage{amssymb}
\usepackage[top=26mm, bottom=26mm, left=32mm, right=32mm]{geometry}
\usepackage[dvips]{color}
\usepackage{bm}
\newenvironment{widequote}{\begin{list}{}
      {\setlength{\rightmargin}{0mm}\setlength{\leftmargin}{5mm}}
      \item[]}{\end{list}}

\newcommand{\vct}[1]{\bm{#1}}
\newcommand{\mtx}[1]{\mathsf{#1}}

\numberwithin{equation}{section}
\numberwithin{figure}{section}
\theoremstyle{definition}
\newtheorem{remark}{Remark}
\numberwithin{remark}{section}
\newtheorem{definition}{Definition}
\numberwithin{definition}{section}

\newcommand{\lsp}{\vspace{3mm}}

\begin{document}

\begin{center}
\textbf{\large A fast solver for Poisson problems on infinite regular lattices}

\lsp

{\small A. Gillman$^1$, P.G. Martinsson$^2$\\
$1$ Department of Mathematics, Dartmouth College, $2$ Department of Applied Mathematics, University of Colorado at Boulder}

\lsp
\begin{minipage}{135mm}
\noindent\textbf{Abstract:}
 The Fast Multipole Method (FMM) provides a highly efficient computational tool for solving constant coefficient partial differential equations (e.g.~the Poisson equation) on infinite domains. The solution to such an equation is given as the convolution between a fundamental solution and the given data function, and the FMM is used to rapidly evaluate the sum resulting upon discretization of the integral.  This paper describes an analogous procedure for rapidly solving elliptic \textit{difference} equations on infinite lattices. In particular, a fast summation technique for a discrete equivalent of the continuum fundamental solution is constructed. 
The asymptotic complexity of the proposed method is $O(N_{\rm source})$, where $N_{\rm source}$ is the number of points subject to body loads.  This is in contrast to FFT 
based methods which solve a lattice Poisson problem at a cost $O(N_{\Omega}\log N_{\Omega})$ independent of $N_{\rm source}$, where $\Omega$ is an artificial rectangular box containing the loaded points and 
$N_{\Omega}$ is the number of points in $\Omega$.  
\end{minipage}
\end{center}

\section{Introduction}
This paper describes an efficient technique for solving Poisson problems
defined on the integer lattice $\mathbb{Z}^{2}$. For simplicity of presentation,
we limit our attention to the equation
\begin{equation}
\label{eq:freespace_lfmm}
[\mtx{A}u](\vct{m}) = f(\vct{m}),\qquad \vct{m} \in \mathbb{Z}^{2},
\end{equation}
where $f = f(\vct{m})$ and $u = u(\vct{m})$ are scalar valued functions
on $\mathbb{Z}^{2}$, and where $\mtx{A}$ is the
so-called \textit{discrete Laplace operator}
\begin{equation}
\label{eq:def_dL}
[\mtx{A}\,u](\vct{m}) = 4u(\vct{m}) - u(\vct{m}+\vct{e}_{1}) - u(\vct{m}-\vct{e}_{1})
                                    - u(\vct{m}+\vct{e}_{2}) - u(\vct{m}-\vct{e}_{2}),
\qquad \vct{m} \in \mathbb{Z}^{2}.
\end{equation}
In (\ref{eq:def_dL}), $\vct{e}_{1} = [1,0]$ and $\vct{e}_{2} = [0,1]$ are the canonical
basis vectors in $\mathbb{Z}^{2}$.  If $f \in L^{1}(\mathbb{Z}^{2})$ and $\sum_{\vct{m} \in \mathbb{Z}^{2}} |f(\vct{m})|  = 0$,
equation (\ref{eq:freespace_lfmm}) is well-posed when coupled with a suitable decay condition for $u$,
see \cite{mdiss} for details.

We are primarily interested in the
situation where the given function $f$ (the \textit{source}) is supported at a
finite number of points which we refer to as \textit{source locations},
and where the function $u$ (the \textit{potential}) is sought at a finite number
of points called \textit{target locations}. While the solution technique is
described for the equation (\ref{eq:freespace_lfmm}) involving the specific operator
(\ref{eq:def_dL}), it may readily be extended
to a broad range of lattice equations involving constant coefficient elliptic
difference operators.

Variations of the equation (\ref{eq:freespace_lfmm}) are perhaps best known as a set of equations
associated with the discretization of elliptic partial differential equations. However, such
equations also emerge in their own right as natural models in a broad range of applications:
random walks \cite{random_walk_cubic}, analyzing the Ising model (in determining
vibration modes of crystals), and many others in engineering mechanics including
 micro-structural models, macroscopic models, simulating fractures \cite{1996_frac,1989_frac} and
as models of periodic truss and frame structures \cite{2001_deshpande_octet_truss,2008_macro,mdiss,2001_gibson_wallach}.

Of particular interest in many of these applications is the
situation where the lattice involves local deviations from perfect periodicity due
to either broken links, or lattice inclusions. The fast technique described in this
paper can readily be modified to handle such situations, see Section \ref{sec:extend_inclusions}.
It may also be modified to handle equations defined on finite subsets of $\mathbb{Z}^{2}$,
with appropriate conditions (Dirichlet / Neumann / periodic) prescribed on the boundary,
see Section \ref{sec:extend_finite_lattices} and \cite{gillman}.

The technique described is a descendant of the Fast Multipole Method (FMM) \cite{rokhlin1987,greengard_thesis,rokhlin1997},
and, more specifically, of ``kernel independent'' FMMs \cite{gimbutas2002,hudson,zorin2004}.
A key application of the original FMM was to rapidly solve the Poisson equation
\begin{equation}
\label{eq:cont_poisson_lfmm}
-\Delta u(\vct{x}) = f(\vct{x}),\qquad \vct{x} \in \mathbb{R}^{2},
\end{equation}
which is the continuum analog of (\ref{eq:freespace_lfmm}). The FMM exploits the fact that
the analytic solution to (\ref{eq:cont_poisson_lfmm}) takes the form of a convolution
\begin{equation}
\label{eq:cont_soln}
u(\vct{x}) = \int_{\mathbb{R}^{2}}\phi_{\rm cont}(\vct{x}-\vct{y})\,f(\vct{y})\,d\vct{y},
\end{equation}
where $\phi_{\rm cont}$ is the fundamental solution of the Laplace operator,
\begin{equation}
\label{eq:def_phi_cont}
\phi_{\rm cont}(\vct{x}) = -\frac{1}{2\pi}\log|\vct{x}|.
\end{equation}
If the source function $f$ corresponds to a number of point charges
$\{q_{j}\}_{j=1}^{N}$ placed at locations $\{\vct{x}_{j}\}_{j=1}^{N}$,
and if the potential $u$ is sought at same set of locations,
then the convolution (\ref{eq:cont_soln}) simplifies to the sum
\begin{equation}
\label{eq:cont_sum}
u_{i} = \underset{j\neq i}{\sum_{j = 1}^N} \phi_{\rm cont}(\vct{x}_{i} - \vct{x}_{j})\,q_{j},
\qquad i = 1,\,2,\,\dots,\,N.
\end{equation}
While direct evaluation of (\ref{eq:cont_sum}) requires $O(N^{2})$ operations since
the kernel is dense, the FMM constructs an approximation to the potentials
$\{u_{i}\}_{i=1}^{N}$
in $O(N)$ operations. Any requested approximation error $\varepsilon$ can be
attained, with the constant of proportionality in the $O(N)$ estimate depending
only logarithmically on $\varepsilon$.

In the same way that the FMM can be said to rely on the fact that the Poisson
equation (\ref{eq:cont_poisson_lfmm}) has the explicit analytic solution (\ref{eq:cont_soln}),
the techniques described in this paper can be said to rely on the fact that
the lattice Poisson equation (\ref{eq:freespace_lfmm}) has an explicit analytic solution
in the form
\begin{equation}
\label{eq:basic_lfmm}
u(\vct{m}) = [\phi * f](\vct{m}) = \sum_{\vct{n} \in \mathbb{Z}^{2}}\phi(\vct{m}-\vct{n})\,f(\vct{n}).
\end{equation}
where $\phi$ is a fundamental solution for the discrete Laplace operator
(\ref{eq:def_dL}). This fundamental solution is known analytically
\cite{duffin1953,mdiss,m2001_LGF_exp,gillman}
via the normalized Fourier integral
\begin{equation}
\label{eq:defphi_lfmm}
\phi(\vct{m}) =
\frac{1}{(2\,\pi)^{2}}\int_{-\pi}^{\pi}\int_{-\pi}^{\pi}\frac{\cos(t_{1}m_{1} + t_{2}m_{2}) - 1}
{4\,\sin^2(t_1/2)+4\,\sin^2(t_2/2)}\,dt_{1}\,dt_{2},\qquad \vct{m} = [m_{1},\,m_{2}]\in \mathbb{Z}^{2}.
\end{equation}

This paper presents an adaptation of the original Fast Multipole Method that enables it
to handle discrete kernels such as (\ref{eq:defphi_lfmm}) and to exploit
 accelerations that are possible due the geometric restrictions
present in the lattice case.  The method extends directly to
any problem that can be solved via convolution with a discrete fundamental solution.
The technique for numerically evaluating (\ref{eq:defphi_lfmm}) extends directly to
other kernels, see Section \ref{sec:eval}.

While we are not aware of any previously published techniques for rapidly
solving the free space problem (\ref{eq:freespace_lfmm}) (or, equivalently, for
evaluating  (\ref{eq:basic_lfmm})), there exist very fast solvers for the closely
related case of lattice Poisson equations defined on rectangular subsets of
$\mathbb{Z}^{2}$ with periodic boundary conditions. Such equations become
diagonal when transformed to Fourier  space, and may consequently be solved
very rapidly via the FFT. The computational time $T_{\rm fft}$ required by
such a method satisfies
\begin{equation}
\label{eq:Tfft}
T_{\rm fft} \sim N_{\rm domain}\log N_{\rm domain}
\qquad\mbox{as}\ N_{\rm domain} \rightarrow \infty,
\end{equation}
where $N_{\rm domain}$ denotes the number of lattice nodes in the smallest rectangular
domain holding all source locations, and where the constant of proportionality is very small. Similar
complexity, sometimes without the logarithmic factor, and with fewer restrictions
on the boundary conditions, may also be achieved via multigrid methods \cite{mccormick_multigrid}.

The principal contribution of the present work is that the computational time
$T_{\rm FMM}$ required by the method described here has asymptotic complexity
\begin{equation}
\label{eq:Tfmm}
T_{\rm FMM} \sim N_{\rm sources},
\qquad\mbox{as}\ N_{\rm sources}\rightarrow \infty,
\end{equation}
where $N_{\rm sources}$ denotes the number of lattice nodes that are loaded
(assuming that the solution is sought only at the source points).
In a situation where the source points
are relatively densely distributed in a rectangle, we would have
$N_{\rm domain} \approx N_{\rm sources}$
and there would be no point in using the new method (in fact, an FFT based
method is in this case significantly faster since the constant of
proportionality in (\ref{eq:Tfft}) is smaller than that in (\ref{eq:Tfmm})).
However, if the source and target points are relatively sparsely distributed
in the lattice, then the estimate (\ref{eq:Tfmm}) of the new method is
clearly superior to that of (\ref{eq:Tfft}) for an FFT based method. As
demonstrated in Section \ref{sec:numerics_lfmm}, very significant gains in
speed can be achieved. Perhaps even more importantly, much larger problems
can be handled since an FFT based method requires that the potential on
all $N_{\rm domain}$ nodes be held in memory.

\begin{figure}[h]
\centering
\includegraphics[width=0.45\textwidth]{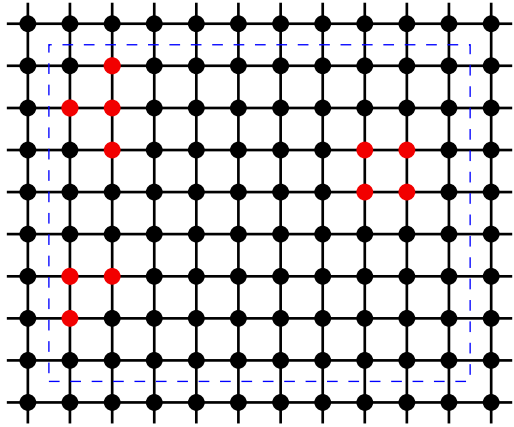}
\caption{A subset of the infinite lattice $\mathbb{Z}^{2}$. The $N_{\rm sources}=11$
red nodes are loaded. The smallest rectangle holding all sources is marked with a
dashed blue line. It has $N_{\rm domain}=80$ nodes.}
\label{fig:lattice}
\end{figure}

\begin{widequote}
\textit{Example:}  The distinction between $N_{\rm domain}$ in (\ref{eq:Tfft})
and $N_{\rm sources}$ in (\ref{eq:Tfmm}) can be illustrated with the toy example
shown in Figure \ref{fig:lattice}. The figure illustrates a portion of
an infinite lattice in which $N_{\rm source}=11$ nodes have been loaded.
A rectangular domain covering these loads is marked with a blue dashed line
and holds $N_{\rm domain}=80$ nodes. Clearly $N_{\rm sources} = 11 \ll 80 = N_{\rm domain}$.
A solution strategy for (\ref{eq:freespace_lfmm}) based on the FFT or multigrid
would involve all $N_{\rm domain}$ nodes inside the rectangle. In contrast,
the lattice fundamental solution allows the solution task to be reduced to evaluating the sum
(\ref{eq:basic_lfmm}) which involves an $N_{\rm sources} \times N_{\rm sources}$
dense coefficient matrix.
\end{widequote}

\section{Review of fast summation techniques}
\label{sec:fmmreview}

In this section, we briefly outline the basic ideas behind the Fast Multipole Method,
and then describe the modifications required to evaluate a lattice sum such as
(\ref{eq:basic_lfmm}).
Our presentation assumes some familiarity with Fast Multipole Methods; for an introduction,
see, \textit{e.g.}, \cite{leslie_FMM_notes,greengard_thesis}.
As a model problem, we consider the task of evaluating the sum
\begin{equation}
\label{eq:finite_sum}
u_{i} = \sum_{j=1}^{N} \phi(\vct{x}_{i}\,-\,\vct{x}_{j})\,q_{j},
\end{equation}
where $\{\vct{x}_{i}\}_{i=1}^{N}$ is a set of $N$ points in the plane,
where $\{q_{i}\}_{i=1}^{N}$ is a set of $N$ given real numbers called \textit{charges},
where $\{q_{i}\}_{i=1}^{N}$ is a set of $N$ sought real numbers called \textit{potentials},
and where $\phi\,\colon\,\mathbb{R}^{2}\times\mathbb{R}^{2} \rightarrow \mathbb{R}$ is a \textit{kernel function.}

\begin{figure}
\centering
\setlength{\unitlength}{1mm}
\begin{picture}(55,55)
\put(00,00){\includegraphics[width=55mm]{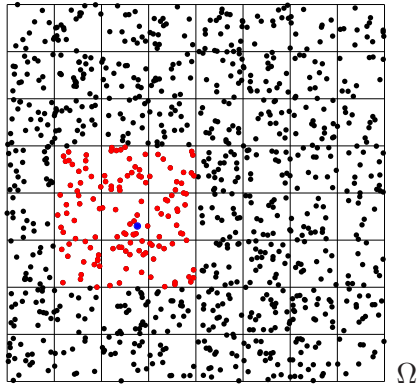}}
\put(54,01){$\Omega$}
\end{picture}
\caption{Geometry of the $N$-body problem in Section \ref{sec:fmmreview}. Source $i$ is blue,
the sources in $J_{i}^{\rm near}$ as defined by (\ref{eq:def_Jnear}) are red.}
\label{fig:fmm_cartoon}
\end{figure}

For simplicity, we consider in this review only the case where the sources are more or less uniformly
distributed in a computational box $\Omega$ in the sense that $\Omega$ can be split into
equi-sized small boxes, called \textit{leaves}, in such a way that each small box holds
about the same number of sources. We let $N_{\rm leaf}$ denote an upper bound for the
number of sources held in any leaf. Then the sum (\ref{eq:finite_sum}) can be split into
two parts
$$
u_{i} = u^{\rm near}_{i} + u^{\rm far}_{i},
$$
where the \textit{near-field} is defined by
\begin{equation}
\label{eq:u_near_lfmm}
u^{\rm near}_{i} = \sum_{j \in J_{i}^{\rm near}\backslash\{i\}}
\phi(\vct{x}_{i} - \vct{x}_{j})\,q_{j},
\qquad i = 1,\,2,\,\dots,\,N,
\end{equation}
where $J_{\rm i}^{\rm near}$ is an index list marking all sources that lie either in the same box as
charge $i$, or in a box that is directly adjacent to the box holding source $i$,
\begin{equation}
\label{eq:def_Jnear}
J_{i}^{\rm near} = \{j\,\colon\,\vct{x}_{i}\mbox{ and }\vct{x}_{j}\mbox{ are located in the same leaf or in directly adjacent leaves}\}.
\end{equation}
The definition of $J_{i}^{\rm near}$ is illustrated in Figure \ref{fig:fmm_cartoon}.
The \textit{far-field} is then defined by
\begin{equation}
\label{eq:u_far_lfmm}
u^{\rm far}_{i} = \sum_{j \notin J_{i}^{\rm near}}
\phi(\vct{x}_{i} - \vct{x}_{j})\,q_{j},
\qquad i = 1,\,2,\,\dots,\,N,
\end{equation}

The \underline{near-field} (\ref{eq:u_near_lfmm}) can now be directly evaluated at low cost
since at most $9\,N_{\rm leaf}$ sources are near any given source. In the
lattice case, this step could potentially be rendered expensive by the fact that the
kernel is known only via the Fourier integral (\ref{eq:defphi_lfmm}) which is quite costly
to evaluate via quadrature. We describe in
Section \ref{sec:eval} how this step may be accelerated by pre-computing and storing
the values of $\phi(\vct{m})$ for all small values of $\vct{m}$ and then using an
asymptotic expansion for large $\vct{m}$. We observe that the local evaluation
gets particularly effective whenever the number of lattice cells along the side of any
leaf box is bounded by some fixed number $L$ of moderate size (say $L \leq 1000$).
In this case, there is in the lattice situation only $16\,L^2$ possible
relative positions of two charges that are near each other which means that
evaluation of the kernel for the near-field calculations amounts to simply a table lookup.
(In fact, due to symmetries, only $2\,L^{2}$ values need to be stored.)

The \underline{far-field} (\ref{eq:u_far_lfmm}) is as in the classical FMM evaluated
via the computation of so-called \textit{multipole expansions} and \textit{incoming expansions}.
These in turn are constructed via a hierarchical procedure on a quad-tree such
as the one shown in Figure \ref{fig:binary_tree_numbered_2D}. With the development
of so-called \textit{kernel independent FMMs}, the multipole expansions of the original
FMM were superseded by more general representations valid for a broad range of kernels.
The bulk of this paper consists of a description of such a kernel independent FMM,
adapted to exploit geometrical restrictions imposed in the lattice case.
Section \ref{sec:expansions} reviews a technique for compactly representing
charges and potentials, and Section \ref{sec:construction} describes how it can be adapted to the
particular case of lattice equations.
Section \ref{sec:tree} introduces notation for handling quad-trees,
Section \ref{sec:operators} describes the so-called \textit{translation operators},
then the full lattice FMM is described in Section \ref{sec:algorithm}.

\section{Evaluation of the lattice fundamental solution}
\label{sec:eval}

The numerical evaluation of the function $\phi$ in (\ref{eq:defphi_lfmm}) requires some care since the integrand has a singularity
at the origin and gets highly oscillatory when $|\vct{m}|$ is large. The latter issue can be handled quite easily since a highly
accurate asymptotic expansion of $\phi(\vct{m})$ as $|\vct{m}|\rightarrow\infty$ is known, see Section \ref{sec:m_large}.
When $|\vct{m}|$ is small, quadrature and Richardson extrapolation may be used to compute $\phi(\vct{m})$ to very high
accuracy, see Section \ref{sec:m_small}.
We note that in the regime where $|\vct{m}|$ is small, there is only a limited
number of possible values of $\vct{m}$, and the corresponding values of $\phi(\vct{m})$
can be pre-computed and stored. Consequently, evaluating $\phi(\vct{m})$ in the near-field
amounts simply to a table look-up, which is very fast.

\subsection{Evaluation of fundamental solution for $|\vct{m}|$ large}
\label{sec:m_large}
It has been established (see \textit{e.g.}~\cite{duffin1953,duffin1958,maradudin1960,mdiss,m2001_LGF_exp})
that as $|\vct{m}| \rightarrow \infty$, the fundamental solution $\phi$ defined by (\ref{eq:defphi_lfmm}) has the asymptotic expansion
\begin{multline}
\label{eq:phi_asympt_lfmm}
\phi(\vct{m}) = -\frac{1}{2\pi}\left(\log|\vct{m}|+\gamma + \frac{\log  8}{2}\right) +
\frac{1}{24\pi}\frac{m_{1}^{4} - 6m_{1}^{2}m_{2}^{2} + m_{2}^{4}}{|\vct{m}|^{6}}\\
+ \frac{1}{480\pi}\frac{43 m_{1}^{8} - 772 m_{1}^{6}m_{2}^{2} +
1570m_{1}^{4}m_{2}^{4} - 772 m_{1}^{2}m_{2}^{6} + 43 m_{2}^{8}}{|\vct{m}|^{12}} +
O(1/|\vct{m}|^{6}).
\end{multline}
The number $\gamma$ is the Euler constant ($\gamma = 0.577206\cdots$).

For $|\vct{m}|$ large, we approximate $\phi$ by dropping the $O(1/|\vct{m}|^{6})$
term off the asymptotic expansion.  We found that for $|\vct{m}|>30$ the expansion
(\ref{eq:phi_asympt_lfmm}) is accurate to at least $10^{-12}$.

The asymptotic expansion (\ref{eq:phi_asympt_lfmm}) is valid for the simple square
lattice only. However, there is a simple process for constructing analogous
expansions for fundamental solutions associated with a very broad class of
constant coefficient elliptic difference operators \cite{m2001_LGF_exp}. The process can be automated
and executed using symbolic software such as Maple \cite{mdiss}.

\subsection{Evaluation of fundamental solution for $|\vct{m}|$ small}
\label{sec:m_small}

When $|\vct{m}|$ is small enough that the asymptotic expansion provides insufficient
accuracy, we approximate the integral (\ref{eq:defphi_lfmm}) using a two-step quadrature procedure:
First, the domain $[-\pi,\,\pi]^{2}$ is split into $n\times n$ equisized boxes
where $n$ is an odd number chosen so that each box holds about one oscillation of the integrand
(in other words, $n \approx |\vct{m}|$). For each box not containing the origin, the integral is approximated
using a Cartesian Gaussian quadrature with $20 \times 20$ nodes. This leaves us with the task
of evaluating the integral
$$
g(a) = \frac{1}{(2\,\pi)^{2}}\int_{-a}^{a}\int_{-a}^{a}\frac{\cos(t_{1}m_{1} + t_{2}m_{2}) - 1}
{4\,\sin^2(t_1/2)+4\,\sin^2(t_2/2)}\,dt_{1}\,dt_{2},
$$
where $a = \pi/n$ denotes the size of the center box. Now observe that
\begin{equation}
\label{eq:Richardson_formula}
g(a) =
\sum_{n=0}^{\infty}\left(g\left(\frac{a}{2^{n}}\right) - g\left(\frac{a}{2^{n+1}}\right)\right) =
\sum_{n=0}^{\infty}\frac{1}{(2\,\pi)^{2}}\int_{\Omega_{n}}\frac{\cos(t_{1}m_{1} + t_{2}m_{2}) - 1}
{4\,\sin^2(t_1/2)+4\,\sin^2(t_2/2)}\,dA,
\end{equation}
where
$$
\Omega_{n} = [2^{-n}\,a,\,2^{-n}\,a]^{2} \backslash [2^{-n-1}\,a,\,2^{-n-1}\,a]^{2},\qquad n = 1,\,2,\,3,\,\dots
$$
is a sequence of annular domains whose union is the square $[-a,\,a]^{2}$.
All integrals in (\ref{eq:Richardson_formula}) involve non-singular integrands,
and can easily be evaluated via Gaussian quadratures. (We split each
$\Omega_{n}$ into eight rectangular regions and use a $20\times 20$ point Gaussian
quadrature on each.) Using
Richardson extrapolation to accelerate the convergence, it turns out that only about 14 terms
are needed to evaluate the sum (\ref{eq:Richardson_formula}) to a precision of $10^{-14}$.

\begin{remark}
The particular integral (\ref{eq:defphi_lfmm}) can be evaluated via a short-cut since it is possible to evaluate the
integral over $t_{1}$ analytically, and then use quadrature
only for the resulting (non-singular) integral over $t_{2}$,
see \cite{mdiss}. Similar tricks are likely possible in many situations involving mono-atomic lattices. However, we prefer
to not rely on this approach since it does not readily generalize to vector valued problems (such as those
associated with mechanical lattice problems) or multi-atomic lattices.
\end{remark}

\section{Outgoing and incoming expansions}
\label{sec:expansions}
In this section, we present techniques for efficiently approximating the far-field $u^{\rm far}_{i}$
to any given positive precision $\varepsilon$.
The parameter $\varepsilon$ can be tuned to balance the
computational cost verses the accuracy. In the numerical examples reported in Section \ref{sec:numerics_lfmm},
$\varepsilon = 10^{-10}$.

\subsection{Interaction ranks}
\label{sec:mpole}
\begin{figure}[ht]
\centering
\setlength{\unitlength}{1mm}
\begin{picture}(50,50)
\put(-15,0){\includegraphics[height=50mm]{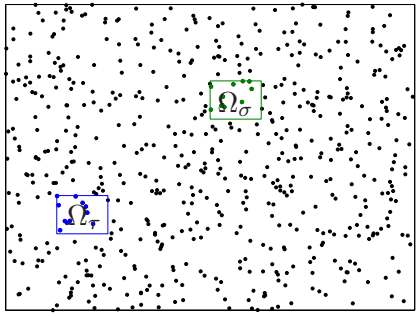}}
\put(1,17){ $\Omega_{\tau}$}
\put(21,32){ $\Omega_{\sigma}$}
\end{picture}
\caption{\label{fig:source} Illustration of source box $\Omega_\tau$ and target box $\Omega_\sigma$.}
\end{figure}

An essential component of the classical FMM
is an efficient technique for representing potentials and source
distributions via ``expansions'' of different kinds. To illustrate
the concept, let us consider a simplified problem in which a number of
sources are placed in a ``source box'' $\Omega_{\tau}$, and the
potential induced by these sources is to be evaluated at a number
of locations in a ``target box'' $\Omega_{\rm \sigma}$. The
orientation of the boxes is shown in Figure \ref{fig:source}. To be precise, we
suppose that sources $\{q_{j}^{\tau}\}_{j=1}^{N}$ are placed at
locations $\{\vct{x}_{j}^{\tau}\}_{j=1}^{N} \subset \Omega_{\rm \tau}$,
and that we seek the potentials $\{u_{i}^{\sigma}\}_{i=1}^{M}$ induced at some
locations $\{\vct{x}_{i}^{\sigma}\}_{i=1}^{M} \subset \Omega_{\rm \sigma}$,
\begin{equation}
\label{eq:twobox}
u_{i}^{\sigma} = \sum_{j=1}^{N}\Phi(\vct{x}_{i}^{\sigma},\,\vct{x}_{j}^{\tau})\,q_{j}^{\tau},
\qquad i = 1,\,2,\,\dots,\,M.
\end{equation}
In this review of the classical FMM, the kernel $\Phi$ is defined by
$$
\Phi(\vct{x},\vct{y}) = \phi_{\rm cont}(\vct{x} - \vct{y}) = -\frac{1}{2\pi}\log|\vct{x} - \vct{y}|,
$$
where $\phi_{\rm cont}$ is the
fundamental solution of the Laplace equation.
For convenience, we write (\ref{eq:twobox}) as a matrix-vector product
\begin{equation}
\label{eq:matvec}
\vct{u}^{\sigma} = \mtx{A}^{\sigma,\tau}\,\vct{q}^{\tau},
\end{equation}
where $\vct{u}^{\sigma} = [u_{i}^{\sigma}]_{i=1}^{M}$ and
$\vct{q}^{\tau} = [q_{j}^{\tau}]_{j=1}^{N}$, and where $\mtx{A}^{\sigma,\tau}$ is
the $M\times N$ matrix with entries
$$
\mtx{A}_{ij}^{\sigma,\tau} = \Phi(\vct{x}_{i}^{\sigma},\,\vct{x}_{j}^{\tau}).
$$

A key observation underlying the FMM is that to any finite precision $\varepsilon$,
the rank of a matrix such as $\mtx{A}^{\sigma,\tau}$ is bounded independently of the numbers $M$
and $N$ of targets and sources in the two boxes. In fact, the $\varepsilon$-rank $P$
of $\mtx{A}^{\sigma,\tau}$ satisfies
$$
P \lesssim \log(1/\varepsilon),\quad\mbox{as}\ \varepsilon \rightarrow 0.
$$
The constant of proportionality depends on the geometry of the boxes,
but is typically very modest.
As a consequence of this rank deficiency, it is possible to factor the matrix
$\mtx{A}^{\sigma,\tau}$, say
\begin{equation}
\label{eq:general_factorization}
\begin{array}{ccccccc}
\mtx{A}^{\sigma,\tau} & \approx & \mtx{B} & \mtx{C},\\
M\times N             &         & M\times P        & P \times N
\end{array}
\end{equation}
and then to evaluate the potential $\vct{u}^{\tau}$ in two steps:
\begin{equation}
\label{eq:fastfactor}
\vct{v} = \mtx{C}\,\vct{q}^{\tau},
\qquad
\vct{u} \approx \mtx{B}\,\vct{v}.
\end{equation}
The cost of evaluating $\vct{u}$ via (\ref{eq:fastfactor}) is
$O((M+N)\,P)$, which should be compared to the $O(M\,N)$ cost
of evaluating $\vct{u}$ via (\ref{eq:matvec}).

\subsection{Formal definitions of outgoing and incoming expansions}
In the classical FMM, a ``multipole expansion'' for a box
is a short vector from which the potential caused by all charges in the box
can be evaluated; it can be viewed as a compressed representation of all the
charges inside the box.
In this section, we introduce the ``outgoing expansion'' as a
generalization of this idea that allows representations other than
classical multipole expansions to be incorporated. The
``incoming expansion'' is analogously introduced to generalize the
concept of a ``local expansion.''

\lsp

\begin{figure}
\centering
\setlength{\unitlength}{1mm}
\begin{picture}(100,45)
\put(20,00){\includegraphics[height=45mm]{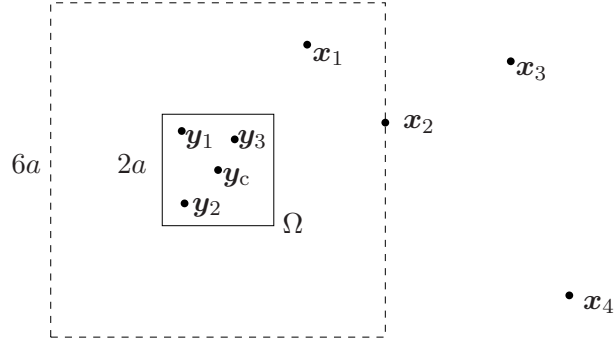}}
\put(15,22){$6a$}
\put(29,22){$2a$}
\put(43,21){$\vct{y}_{\rm c}$}
\put(38,26){$\vct{y}_{1}$}
\put(39,17){$\vct{y}_{2}$}
\put(45,26){$\vct{y}_{3}$}
\put(55,37){$\vct{x}_{1}$}
\put(67,28){$\vct{x}_{2}$}
\put(82,35){$\vct{x}_{3}$}
\put(91,04){$\vct{x}_{4}$}
\put(51,14){$\Omega$}
\end{picture}
\caption{Illustration of \textit{well-separated points}. Any point on or outside of the dashed square
is \textit{well-separated} from $\Omega$. Consequently, the points $\vct{x}_{2}$, $\vct{x}_{3}$, and $\vct{x}_{4}$
are well-separated from $\Omega$, but $\vct{x}_{1}$ is not.}
\label{fig:wellsep}
\end{figure}

\noindent
\underline{\textit{Well-separated boxes:}} Let $\Omega$ be a box with
side length $2a$ and center $\vct{c}$ as shown in Figure \ref{fig:wellsep}.
We say that a point $\vct{x}$ is \textit{well-separated} from $\Omega$ if
it lies outside the square of
side length $6a$ centered at $\vct{c}$. We say that two boxes $\Omega$
and $\Omega'$ are \textit{well-separated} if every point in $\Omega'$ is
well-separated from $\Omega$, and vice versa.

\lsp

\noindent
\underline{\textit{Outgoing expansion:}}
Let $\Omega$ be a box containing a set of sources.
We say that a vector $\hat{\vct{q}}$ is an \textit{outgoing expansion}
for $\Omega$ if the potential caused by the sources in $\Omega$ can
be reconstructed from $\hat{\vct{q}}$ to within precision $\varepsilon$
at any point that is well-separated from $\Omega$.

\lsp

\noindent
\underline{\textit{Incoming expansion:}}
Let $\Omega$ be a box in which a potential has been induced by a set
of sources located at points that are well-separated from $\Omega$.
We say that a vector $\vct{\hat{u}}$ is an \textit{incoming expansion}
for $\Omega$ if $u$ can be reconstructed from $\hat{\vct{u}}$
to within precision $\varepsilon$.

\subsection{Charge basis}
\label{sec:chargebasis}
The cost of computing a factorization such as (\ref{eq:general_factorization})
using a generic linear algebraic technique such as QR is $O(M\,N\,P)$, which would
negate any savings obtained when evaluating the matrix-vector product
(unless a very large number of matrix-vector products involving the same
source and target locations is required). Fortunately, it is possible in many
environments to construct such factorizations much faster.  The classical FMM
uses multipole expansions.  As an alternative, an approach based
on so-called ``proxy charges'' has recently been developed \cite{zorin2004}.
It has been demonstrated \cite{hudson,tygert2005} that for any given box $\Omega$, it is
possible to find a set of locations $\hat{\vct{Y}} = \{\hat{\vct{y}}_{p}\}_{p=1}^{P} \subset \Omega$
with the property that sources placed at these points can to high accuracy replicate
any potential caused by a source distribution in $\Omega$.  The number of points $P$
required is given in Table \ref{tab:ranks}. To be precise, given
any set of points $\vct{Y} = \{\vct{y}_{j}\}_{j=1}^{N}\subset \Omega$ and any sources
$\vct{q} = \{q_{j}\}_{j=1}^{N}$, we can find ``equivalent charges''
$\hat{\vct{q}} =\{\hat{q}_{p}\}_{p=1}^{P}$
such that
\begin{equation}
\label{eq:proxy}
\sum_{j=1}^{N}\phi_{\rm cont}(\vct{x}\,-\,\vct{y}_{j})\,q_{j} \approx
\sum_{p=1}^{P}\phi_{\rm cont}(\vct{x}\,-\,\hat{\vct{y}}_{p})\,\hat{q}_{p},
\end{equation}
whenever $\vct{x}$ is well-separated from $\Omega$. The approximation (\ref{eq:proxy})
holds to some preset (relative) precision $\varepsilon$.
Moreover, the map from $\vct{q}$ to $\hat{\vct{q}}$ is linear, and there exists
a matrix $\mtx{T}_{\rm ofs} = \mtx{T}_{\rm ofs}(\hat{\vct{Y}},\,\vct{Y})$ such that
\begin{equation}
\label{eq:T_ofs}
\hat{\vct{q}} = \mtx{T}_{\rm ofs}\,\vct{q},
\end{equation}
 where ``ofs'' is an abbreviation of ``outgoing [expansion] from sources.''

\begin{table}[ht]
\centering
 \begin{tabular}{|c|c|c|c|c|}
\hline
   & \multicolumn{4}{|c|}{$\epsilon$}\\
\hline
  & &  & & \\
  $l $ & $10^{-6}$& $10^{-8}$ & $10^{-10}$& $10^{-13}$\\
\hline
$1/32$ & $19$ & $27$ & $37$ & $49$\\
\hline
$1/16$ & $19$ & $27$ & $36$ & $49$\\
\hline
$1/4$ & $19$ & $28$ & $37$& $ 49$\\
\hline
$1/2$ & $21$ & $29$ & $37$& $ 51$\\
\hline
 \end{tabular}
 \caption{\label{tab:ranks} The number of points $P$ required to replicate the field to accuracy $\epsilon$ for a box $\Omega$
with side length $l$.}
\end{table}

We say that the points $\{\hat{\vct{y}}_{p}\}_{p=1}^{P}$ form an
\textit{outgoing skeleton} for $\Omega$, and that the vector $\hat{\vct{q}}$
is an outgoing expansion of $\Omega$.

In addition, we can find an \textit{incoming skeleton}
$\hat{\vct{X}} = \{\hat{\vct{x}}_{p}\}_{p=1}^{P} \subset \Omega$
with the property that any incoming potential in $\Omega$ can be
interpolated from its values on the incoming skeleton. To be precise,
suppose that $U = U(\vct{x})$ is a potential caused by sources that
are well-separated from $\Omega$, and that $\vct{X} = \{\vct{x}_{i}\}_{i=1}^{M}$
is an arbitrary set of points in $\Omega$. Then there exists a matrix
$\mtx{T}_{\rm tfi} = \mtx{T}_{\rm tfi}(\vct{X},\hat{\vct{X}})$ (`tfi'' stands for ``targets from incoming [expansion]'')
 such that
$$
\vct{u} = \mtx{T}_{\rm tfi}\,\hat{\vct{u}},
$$
where
$$
\vct{u} = [U(\vct{x}_{i})]_{i=1}^{M},
\qquad\mbox{and}\qquad
\hat{\vct{u}} = [U(\hat{\vct{x}}_{p})]_{p=1}^{P}.
$$

When the kernel $\Phi$ is symmetric in the sense that
$\Phi(\vct{x}\,-\,\vct{y}) = \Phi(\vct{y}\,-\,\vct{x})$ for all $\vct{x}$ and $\vct{y}$,
any outgoing skeleton is also an incoming skeleton,
$$
\hat{\vct{X}} = \hat{\vct{Y}}.
$$
Moreover, if the target points equal the source points so that $\vct{X} = \vct{Y}$,
then
$$
\mtx{T}_{\rm tfi} = \left(\mtx{T}_{\rm ofs}\right)^{*}.
$$

Applied to the situation described in Section \ref{sec:mpole},
where a set of sources were placed in a source box $\Omega_{\tau}$, and we sought
to evaluate the potential induced at a set of target points in a box $\Omega_{\sigma}$,
the claims of this section can be summarized by saying that $\mtx{A}^{\sigma,\tau}$ admits
an approximate factorization
$$
\begin{array}{ccccccc}
\mtx{A}^{\sigma,\tau} & \approx & \mtx{T}_{\rm tfi}^{\sigma} &
                                  \mtx{T}_{\rm ifo}^{\sigma,\tau} &
                                  \mtx{T}_{\rm ofs}^{\tau}\\
M\times N             &         & M\times P        & P \times P            & P \times N
\end{array}
$$
where the middle factor is simply a subsampling of the original kernel function
$$
\mtx{T}_{\mbox{ifo},pq}^{\sigma,\tau} = \Phi(\hat{\vct{x}}_{p}^{\sigma},\,\hat{\vct{x}}_{q}^{\tau}).
$$

\begin{remark}
For solving multiple problems involving different source and load distributions that involve the same kernel, one set of skeleton points
may be used for all problems by choosing the skeleton points to lie on the boundary of $\Omega_{\sigma}$ and $\Omega_{\tau}$. The interpolation matrices
$\mtx{T}_{\rm tfi}$ and $\mtx{T}_{\rm ofs}$ need be constructed for each unique set of source and load distributions using the techniques from
\cite{mskel}.  In Section \ref{sec:construction}, we describe this generalization
of the skeletonization process in more detail for the lattice fundamental solution.

\end{remark}

\section{Constructing charge bases for the lattice fundamental solution}
\label{sec:construction}
In this section, we describe how to construct the charge bases for the lattice fundamental solution defined by (\ref{eq:defphi_lfmm}).

From potential theory, we know that to capture the interaction between a set of source points $ \{\vct{m}^\tau_j\}_{j=1}^N$ in box $\Omega_{\tau}$
and all points far from $\Omega_{\tau}$, it is enough to capture the interaction between the source points and a set of ``proxy'' points $F$ that
lie densely on the boundary of a box that is concentric to $\Omega_{\tau}$ and has a boundary that is well-separated from $\Omega_{\tau}$.

We choose the skeleton points to be a subset of the set of all points $Y$ that lie on the boundary of box $\tau$.  Either
rank revealing QR factorization \cite{gu_rrqr} or factorization techniques from \cite{mskel} are applied to the matrix  $\mtx{A}^{F,Y}$ (whose entries
are given  by $\mtx{A}^{F,Y}_{i,j} = \Phi(\vct{m}_i^F\,-\,\vct{m}_j^Y)$)
to determine the rank $P$ and which $P$ points make up the set of skeleton points $\hat{Y}$ of $\Omega_{\tau}$.

Using the skeleton points and the techniques from \cite{mskel}, we find the $P\times N$ matrix $\mtx{T}_{\rm{ofs}}$ such that
\begin{equation} \label{eq:skel} \| \mtx{A}^{F,\tau}-\mtx{A}^{F,\hat{Y}}\mtx{T}_{\rm{ofs}}\| < \epsilon.\end{equation}

We use a similar technique to find the incoming skeleton points and the translation operator $\mtx{T}_{\rm{tfi}}$.

\begin{remark} Because of the smoothness of the kernel, it is not required to use all the points on the boundary of the well-separated box
as ``proxy'' points.  We found it is enough to take 40 points per edge to approximate the far field with accuracy $10^{-10}$.
\end{remark}

\section{Tree structure}
\label{sec:tree}
The separation of variables in the kernel that was described in Section \ref{sec:expansions}
is all that is needed to effectively evaluate a potential field whenever the set of target
locations is well-separated from the set of source locations.
 When the two sets coincide, we need to tessellate the box containing them into smaller boxes, and use the
expansion only for interactions between boxes that are well-separated. In this section, we describe the
simplest such tessellation.

Suppose that we are given a set of points $\{\vct{x}_{i}\}_{i=1}^{N}$
in a box $\Omega$. Given an integer $N_{\rm leaf}$, we pick the smallest
integer $L$ such that when the box $\Omega$ is split into
$4^{L}$ equisized smaller boxes, no box holds more than $N_{\rm leaf}$ points.
These $4^{L}$ equisized small boxes form the
\textit{leaves} of the tree. We merge the leaves by sets of four
to form $4^{L-1}$ boxes of twice the side-length, and then continue
merging by pairs until we recover the original box $\Omega$, which we
call the \textit{root}.

\begin{figure}
\centering
\input{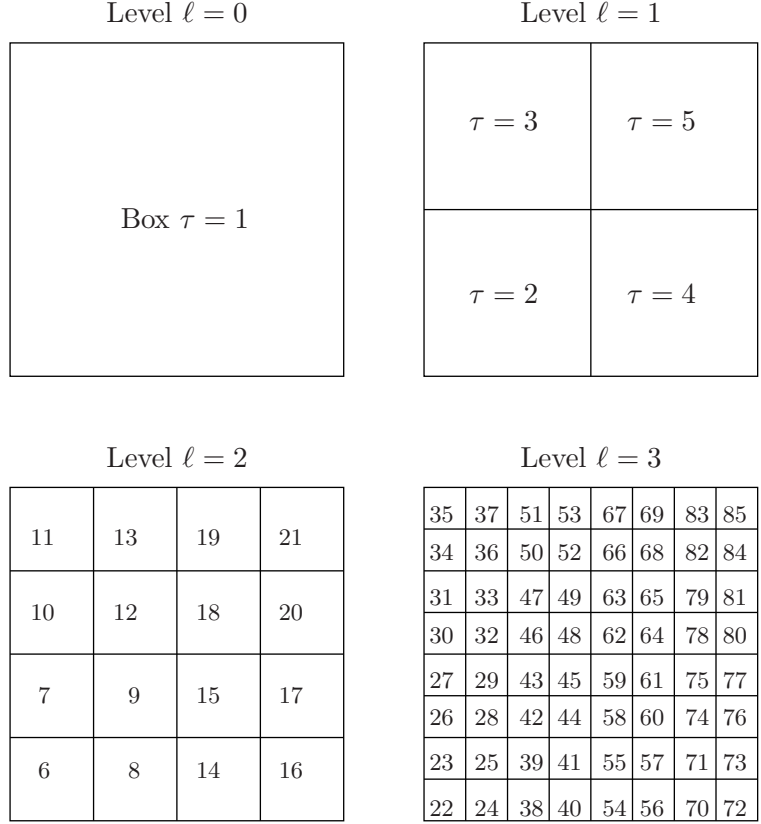}
\caption{A binary tree with 4 levels of uniform refinement. }
\label{fig:binary_tree_numbered_2D}
\end{figure}

The set consisting of all boxes of the same size forms what we
call a \textit{level}. We label the levels using the integer $\ell = 0,\,1,\,2,\,\dots,\,L$,
with $\ell = 0$ denoting the root, and $\ell = L$ denoting the leaves.
See Figure \ref{fig:binary_tree_numbered_2D}.

\begin{definition}
\label{def:tree}
Let $\tau$ be a box in a hierarchical tree.
\begin{itemize}
\item The \underline{parent}
of $\tau$ is the box on the next coarser level that contains $\tau$.
\item The \underline{children}
of $\tau$ is the set $\mathcal{L}_{\tau}^{\rm child}$ of boxes whose parent is $\tau$.
\item The \underline{neighbors} of $\tau$ is the set
$\mathcal{L}^{\rm nei}_{\tau}$ of boxes that
are on the same level as $\tau$ and are directly adjacent to it.
\item The \underline{interaction list} of $\tau$ is the set
$\mathcal{L}_{\tau}^{\rm int}$ of all boxes $\sigma$ such that:
\begin{enumerate}
\item $\sigma$ and $\tau$ are on the same level.
\item $\sigma$ and $\tau$ are not directly adjacent.
\item The parents of $\sigma$ and $\tau$ are directly adjacent.
\end{enumerate}
\end{itemize}
\end{definition}

\lsp

\noindent
\textit{Example:}
For the tree shown in Figure \ref{fig:binary_tree_numbered_2D}, we have, \textit{e.g.},
\begin{align*}
\mathcal{L}_{14}^{\rm child} &= \{54,\,55,\,56,\,57\},\\
\mathcal{L}_{23}^{\rm nei}   &= \{22,\,24,\,25,\,26,\,28\},\\
\mathcal{L}_{59}^{\rm nei}   &= \{36,\,37,\,48,\,58,\,60,\,61,\,70,\,72\},\\
\mathcal{L}_{ 7}^{\rm int}   &= \{11,13,14\colon 21\},\\
\mathcal{L}_{37}^{\rm int}   &= \{22\colon 29,\,30\colon33,\,38\colon 41,\,47,\,49,\,54\colon57,\,60,\,61,\,71,\,72,\,73\}.
\end{align*}

For the moment, we are assuming that the given point distribution
is sufficiently uniform that all the leaves hold roughly the same
number of points. In this case,
$$
L \sim \log\frac{N}{N_{\rm leaf}}.
$$
For non-uniform distributions of points, a uniform subdivision
of $\Omega$ into $4^{L}$ boxes of equal length would be
inefficient since many of the leaves would hold few or no points.
In such cases, adaptive subdivisions should be used \cite{carrier_greengard}.

\section{Translation operators}
\label{sec:operators}
In the FMM, five different
so-called \textit{translation operators} that construct outgoing or incoming expansions are required.
We will, in this section, describe how to construct them, but we
first list which operators we need:

\begin{itemize}

\item[$\mtx{T}_{\rm ofs}^{\tau}$]
\underline{\textit{The outgoing from sources translation operator:}}
Let $\tau$ denote a box holding a set of sources whose values are listed
in the vector $\vct{q}^{\tau}$. The outgoing expansion $\hat{\vct{q}}^{\tau}$
of $\tau$ is then constructed via
$$
\hat{\vct{q}}^{\tau} = \mtx{T}_{\rm ofs}^{\tau}\,\vct{q}^{\tau}.
$$

\lsp

\item[$\mtx{T}_{\rm ofo}^{\tau,\sigma}$]
\underline{\textit{The outgoing from outgoing translation operator:}}
Suppose that a child $\sigma$ of a box $\tau$ holds a source distribution
represented by the outgoing expansion $\hat{\vct{q}}^{\sigma}$.
The far-field caused by these sources can equivalently be
represented by an outgoing representation $\hat{\vct{q}}^{\tau}$ of
the parent, constructed via
$$
\hat{\vct{q}}^{\tau} = \mtx{T}_{\rm ofo}^{\tau,\sigma}\,\hat{\vct{q}}^{\sigma}.
$$

\lsp

\item[$\mtx{T}_{\rm ifo}^{\tau,\sigma}$]
\underline{\textit{The incoming from outgoing translation operator:}}
Suppose that $\tau$ and $\sigma$ are two well-separated boxes, and that $\sigma$
holds a source distribution represented by an outgoing expansion $\hat{\vct{q}}^{\sigma}$.
Then the field in $\tau$ caused by these sources can be represented by an
incoming expansion $\hat{\vct{u}}^{\tau}$ that is constructed via
$$
\hat{\vct{u}}^{\tau} = \mtx{T}_{\rm ifo}^{\tau,\sigma}\,\hat{\vct{q}}^{\sigma}.
$$

\lsp

\item[$\mtx{T}_{\rm ifi}^{\tau,\sigma}$]
\underline{\textit{The incoming from incoming translation operator:}}
Suppose that $\tau$ is the parent of a box $\sigma$.
Suppose further that the incoming expansion $\hat{\vct{u}}^{\tau}$
represents a potential in $\tau$ caused by sources that are all
well-separated from $\tau$. Then these sources are also well-separated
from $\sigma$, and the potential in $\sigma$ can be represented via an
incoming expansion $\hat{\vct{u}}^{\sigma}$ given by
$$
\hat{\vct{u}}^{\sigma} = \mtx{T}_{\rm ifi}^{\tau,\sigma}\,\hat{\vct{u}}^{\sigma}.
$$

\lsp

\item[$\mtx{T}_{\rm tfi}^{\tau}$]
\underline{\textit{The targets from incoming translation operator:}}
Suppose that $\tau$ is a box whose incoming potential is
represented via the incoming representation $\hat{\vct{u}}^{\tau}$.
Then the potential at the actual target points are constructed via
$$
\vct{u}^{\tau} = \mtx{T}_{\rm tfi}^{\tau}\,\hat{\vct{u}}^{\tau}.
$$
\end{itemize}

Techniques for constructing the matrix $T_{\rm ofs}^{\tau}$ were
described in Section \ref{sec:construction}. Since in our case, the kernel is
symmetric (\textit{i.e.}~$\phi(\vct{x}\,-\,\vct{y}) = \phi(\vct{y}\,-\,\vct{x})$ for all
$\vct{x}$ and $\vct{y}$), these techniques immediately give us the targets-from-incoming
translation operator as well, since
$$
\mtx{T}_{\rm tfi}^{\tau} = \left(\mtx{T}_{\rm ofs}^{\tau}\right)^{*}.
$$

We next observe that when charge bases are used, the outgoing-to-incoming
translation operator is simply a sampling of the kernel function,
$$
\mtx{T}_{\mbox{ifo},pq}^{\tau,\sigma} =
\phi(\hat{\vct{x}}_{p}^{\tau}\,-\,\hat{\vct{x}}_{q}^{\sigma}),
\qquad p,\,q = 1,\,2,\,3,\,\dots,\,P,
$$
where
$\{\hat{\vct{x}}_{i}^{\tau}\}_{i=1}^{P}$ and
$\{\hat{\vct{x}}_{j}^{\sigma}\}_{j=1}^{P}$ are the locations of the
skeleton points of $\tau$ and $\sigma$, respectively.

All that remains is to construct $\mtx{T}_{\rm ofo}^{\tau,\sigma}$
and $\mtx{T}_{\rm ifi}^{\sigma,\tau}$. In fact, since the kernel is symmetric,
$$
\mtx{T}_{\rm ifi}^{\sigma,\tau} =
\left(\mtx{T}_{\rm ofo}^{\sigma,\tau}\right)^{*},
$$
and all that actually remains is to construct the matrices $\mtx{T}_{\rm ofo}^{\tau,\sigma}$.
To this end, let $\{\sigma_i\}_{i=1}^{l}$ denote children of box $\tau$.
The construction of $\mtx{T}_{\rm ofo}^{\tau,\sigma_i}$
closely resembles the construction of the $\mtx{T}_{\rm ofs}$ operator described in Section \ref{sec:construction}.
Instead of choosing the skeleton points from the set of all points on the boundary of $\tau$ as was done in the construction
of $\mtx{T}_{\rm ofs}$, we choose the skeleton points for $\tau$ to be a subset of the skeleton points of its children,
 $Y = [\hat{Y}^{\sigma_1},\ldots,\hat{Y}^{\sigma_l}]$.  As in Section \ref{sec:construction}, we define a set of ``proxy'' points
$F$ that are well-separated from $\tau$ and use a factorization technique such as rank revealing QR to determine which points in $Y$
make up the set of skeleton points $\hat{Y}$.  Using the techniques from \cite{mskel}, we find the interpolation matrix $\mtx{S}$ such that
$$
\| \mtx{A}^{P,Y}-\mtx{A}^{P,\hat{Y}}\mtx{S}\| < \epsilon.
$$
The translation operator $\mtx{T}_{\rm{ofo}}^{\tau,\sigma_1}$ is then defined via
$
\mtx{T}_{\rm{ofo}}^{\tau,\sigma_1} = \mtx{S}(:,1:k_1)
$
where $k_1$ is the number of skeleton points of $\sigma_1$, $\mtx{T}_{\rm{ofo}}^{\tau,\sigma_2} = \mtx{S}(:,(k_1+1):(k_1+k_2))$ where $k_2$ is the number of skeleton points of $\sigma_2$, etc.

\section{A lattice Fast Multipole Method}
\label{sec:algorithm}

While the classical FMM derives so-called ``translation operators'' based
on asymptotic expansions of the kernel function, the method we propose
 determines these operators computationally.  In this regard, it is similar to
``kernel independent FMMs'' such as \cite{anderson1992,makino1999,zorin2004}.  Since the kernel is translation invariant,
the computations need be carried out only for a single box on each level.  Thus the
construction of the translation operators is very inexpensive (less than linear complexity).

\subsection{Precomputing skeletons and translation operators}
\label{sec:precompute}

For each level $l$, we define a ``model'' box which is centered at the origin and has the same size as the boxes on level $l$.  The skeleton points
 and the translation operators are found with respect to the model box.

To illustrate the concept, suppose that we are given a source $f$ that is non-zero set of points $\{\vct{m}_{i}\}_{i=1}^{N}$
in a box $\Omega$.  We seek the potential at the source points.

 The pre-computation consist of the following steps:

\begin{enumerate}
 \item Divide $\Omega$ into the tree structure as described in Section \ref{sec:tree}.
\item Construct the lists described in Section \ref{sec:tree}.
\item Construct the skeleton points,  $\mtx{T}_{\rm ofo}$, and $\mtx{T}_{\rm ifi}$ translation operators.  At the lowest level $L$, we construct
the skeleton points for the level $L$ model box using the procedure described in Section \ref{sec:construction}.  For each level $i<L$, we take
four copies of the skeleton points for level $i+1$ shifting them so that each copy makes up one quadrant of the model box for level $i$.  The skeleton
points and the translation operators $\mtx{T}_{\rm{ofo}}$ and $\mtx{T}_{\rm{ifi}}$ are constructed using the technique described in Section
\ref{sec:operators}.
\item Construct the $\mtx{T}_{\rm ifo}$ translation operators.  For each level $i>1$, we construct the $\mtx{T}_{\rm{ifo}}$ translation
operators for the model box.  We assume that the model box is completely surrounded with boxes such that the interaction list has the maximum
number of boxes possible which is 42.  Let $\hat{Y}$ be the outgoing skeleton points and $\hat{X}$ be the incoming skeleton points for the model
 box on level $i$.  For each $j\leq 42$, we shift $\hat{X}$ to be centered at the $j^{\rm{th}}$ possible location for a box on the interaction
list and define
\begin{equation} \tilde{\mtx{T}}^j_{\rm{ifo}} = \mtx{A}^{\hat{X},\hat{Y}}\end{equation}
\end{enumerate}

\begin{remark}
In computing the sum, described in Section \ref{sec:application}, it is easy to use the pre-computed translation operators.
 For example, given a box $\tau$ that has a box $\sigma$ on the interaction list, we identify which $j$ location $\sigma$ is in relative
to $\tau$ and define $\mtx{T}^{\sigma,\tau} = \tilde{\mtx{T}}^j_{\rm{ifo}}$.
\end{remark}

\begin{remark}
For leaf boxes of size less than $8\times 8$
on level $l$, we utilize the fact that there are a finite number of points inside the
 box that are also in $\mathbb{Z}^2$ and construct the translation operator $\mtx{T}^l_{\rm{ofs}}$ for the model box assuming the source
 points are dense. For each box $\tau$ on level $l$ with ${N^\tau}$ sources, we construct an index vector $J^{\tau}$ that notes the
locations of the sources $\{\vct{m}^\tau_j\}_{j=1}^{N^\tau}$ in the dense lattice.  We define
 $\mtx{T}^{\tau}_{\rm{ofs}} = \mtx{T}^l_{\rm{ofs}}(:,J^\tau)$. The translation operator $\mtx{T}^{\tau}_{\rm{tfi}}$ is constructed in a similar manner.
\end{remark}

\subsection{Application}
\label{sec:application}
We have now assembled the tools for computing the sum
(\ref{eq:basic_lfmm}) through two passes through the hierarchical
tree; one upwards, and one downwards.

\lsp

\begin{enumerate}
\item
Sweep over all leaf boxes $\tau$. For each box, construct its
outgoing representation from the values of the sources inside it:
$$
\hat{\vct{q}}^{\tau} = \mtx{T}_{\rm ofs}^{\tau}\,\vct{q}(J^{\tau}).
$$

\lsp

\item
Sweep over all non-leaf boxes $\tau$, going from finer to coarser levels.
Merge the outgoing expansions of the children to construct the outgoing
expansion for $\tau$,
$$
\hat{\vct{q}}^{\tau} = \sum_{\sigma \in \mathcal{L}_{\rm children}^{\tau}}
\mtx{T}_{\rm ofo}^{\tau,\sigma} \hat{\vct{q}}^{\sigma}.
$$

\lsp

\item Loop over all boxes $\tau$. For each box, collect the contributions to its
incoming expansion from boxes in its interaction list:
$$
\hat{\vct{u}}^{\tau} =
\sum_{\sigma \in \mathcal{L}_{\rm int}^{\tau}}
\mtx{T}_{\rm ifo}^{\tau,\sigma}\,\hat{\vct{q}}^{\sigma}.
$$

\lsp

\item Loop over all parent boxes $\tau$, going from coarser levels
to finer. For each box $\tau$, loop
over all children $\sigma$ of $\tau$, and broadcast the
the incoming expansion of $\tau$ to the incoming expansions of
$\sigma$:
$$
\hat{\vct{u}}^{\sigma} = \hat{\vct{u}}^{\sigma} +
\mtx{T}_{\rm ifi}^{\sigma,\tau}\,\hat{\vct{u}}^{\tau}.
$$

\lsp

\item
Sweep over all leaf nodes $\tau$. For each node, form the potential
$\vct{u}^{\tau}$
by evaluating the incoming representation and directly adding the contributions
from the sources inside $\tau$ and in all boxes that are not well-separated
from $\tau$:
$$
\vct{u}^{\tau} = \vct{u}(J^{\tau}) =
\mtx{T}_{\rm tfi}^{\tau}\,\hat{\vct{u}}^{\tau} +
\mtx{A}(J^{\tau},J^{\tau})\,\vct{q}(J^{\tau}) +
\sum_{\sigma \in \mathcal{L}_{\rm nei}^{\tau}}
\mtx{A}(J^{\tau},J^{\sigma})\,\vct{q}(J^{\sigma}).
$$
\end{enumerate}

\subsection{Asymptotic complexity of the proposed scheme}
Since the kernel (\ref{eq:defphi_lfmm}) is separable, the cost of computing the skeleton points and the translation operators
on any level of the quad-tree is $O(P\,M\, |F|)$ where $P$ is the number of skeleton points,
 $M$ is the number of points on the boundary the box, and  $|F|$ is the number of
well-separated proxy nodes used \cite{mskel}.  The cost of solving the least squares problem (\ref{eq:skel}) to find the matrix
$T_{\rm ofs}$ for a leaf box is $O(P^2\,|F| + N\,P\,|F|)$  where $N$ is the number of loaded points in the box.  Hence,
the total complexity of the lattice FMM is $O(N_{\rm source})$.

Also notice that the memory needed to store the precomputed information is $O(N_{\rm source})$.

\section{Numerical examples}
\label{sec:numerics_lfmm}

In this section, we show that the lattice FMM speed compares favorably
to FFT based techniques except for situations where the
source points populate the majority of some rectangular computational box.
We also show that the amount of memory required scales linearly with the number of source terms.

All experiments are run on a Dell desktop computer with 2GB of RAM and an Intel Pentium 4
3.4GHz dual processor. The method was run at a requested relative precision of $10^{-10}$.
The techniques were implemented rather crudely in Matlab, which means that significant further
gains in speed should be achievable.

We consider the lattice Poisson problem
\begin{equation}
\label{eq:basic1}
[\mtx{A}u](\vct{m}) = f(\vct{m}),
\end{equation}
where the points where $f(\vct{m})$ is non-zero are confined to an $n \times n$ square subdomain $\Omega $ of $\mathbb{Z}^2$.
The FFT produces a slightly different solution than the lattice FMM since it enforces periodic boundary conditions, but this is
not important for our purposes.  We suppose throughout that $n$ is a power of two to make the comparison as favorable to the FFT as possible. We
let $T_{\rm fft}$ denote the time required by the FFT, and $T_{\rm FMM}$ the time for the FMM.

\begin{figure}[ht]
\centering
\setlength{\unitlength}{1mm}
\begin{picture}(150,56)
\put(00,5){\includegraphics[width=40mm]{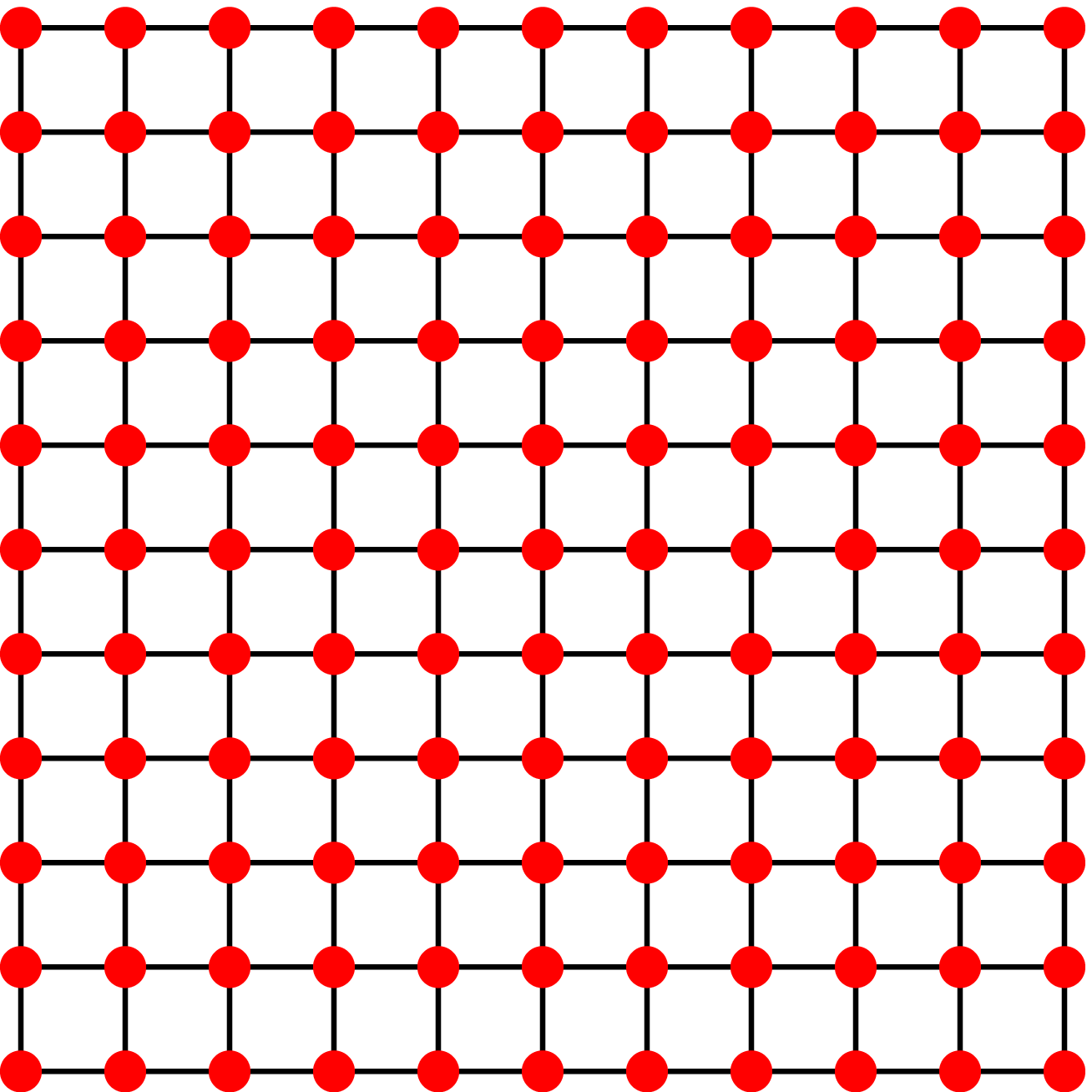}}
\put(50,5){\includegraphics[width=40mm]{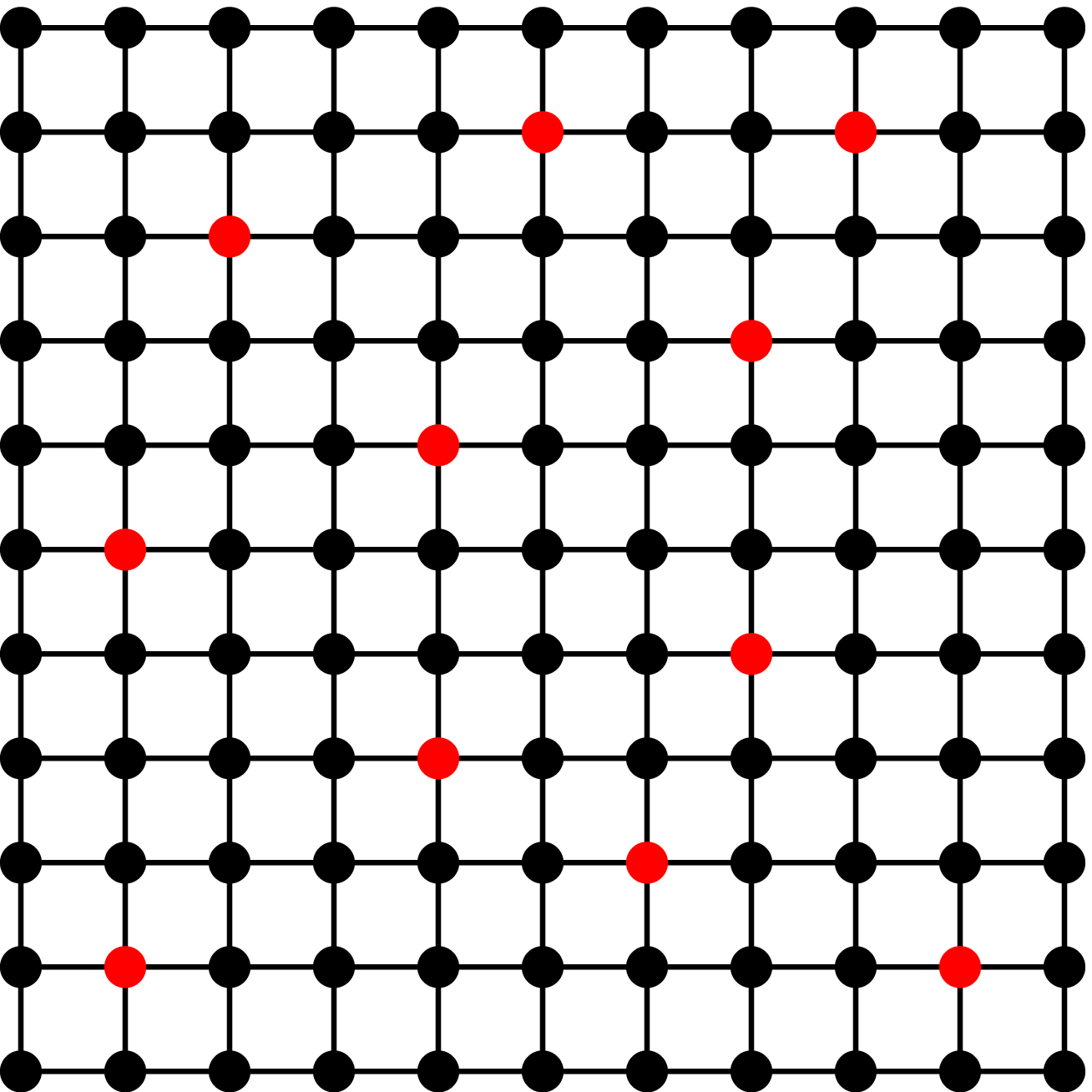}}
\put(100,5){\includegraphics[width=40mm]{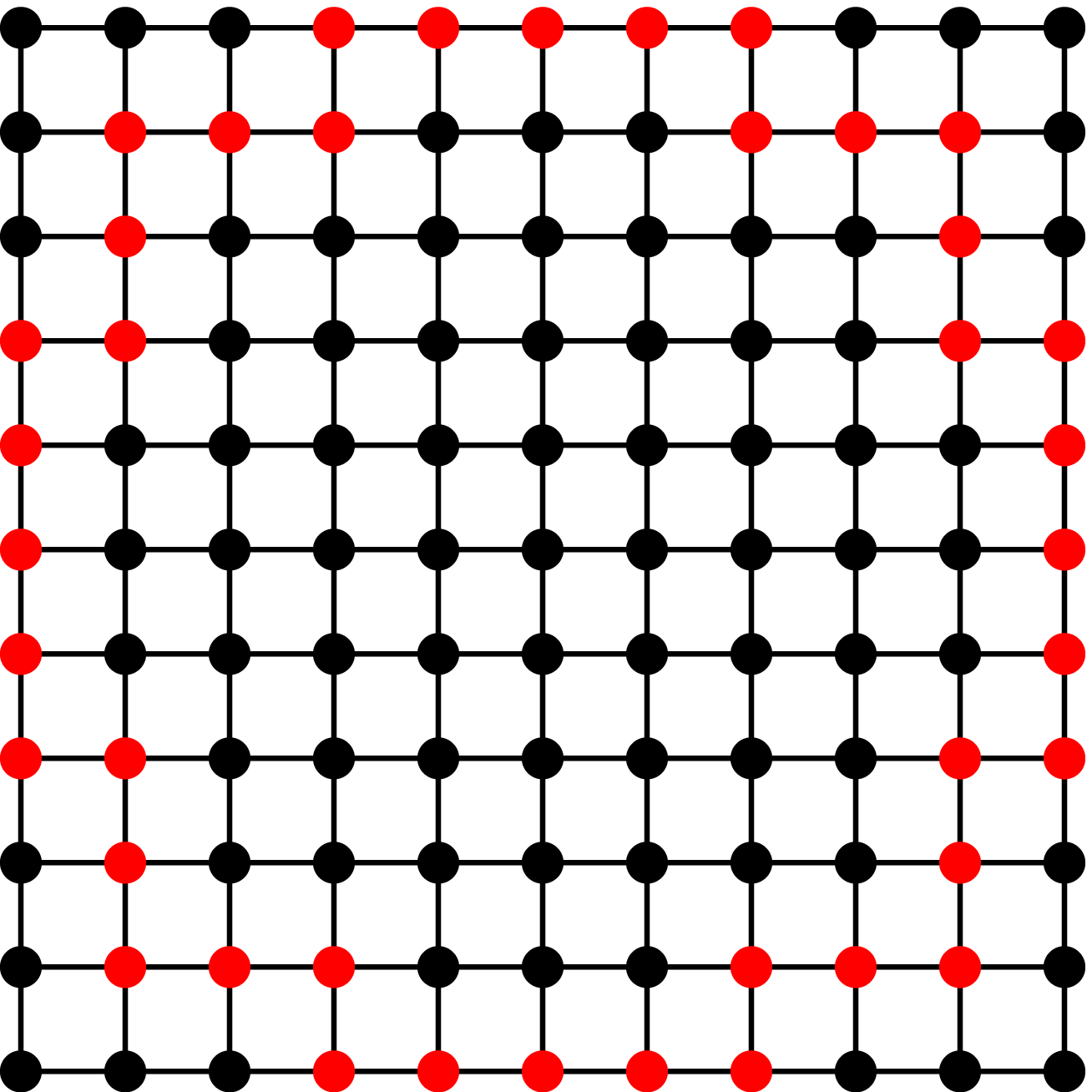}}
\put(10,00){(a) Dense}
\put(60,00){(b) Random}
\put(105,00){(c) Embedded Circle}
\end{picture}
\caption{\label{fig:geometries} Illustration of load distributions for the three experiments.  Red dots are the source points.}
\end{figure}

\begin{figure}[ht]
\centering
\setlength{\unitlength}{1mm}
\begin{picture}(160,85)
\put(-5,5){\includegraphics[width=85mm]{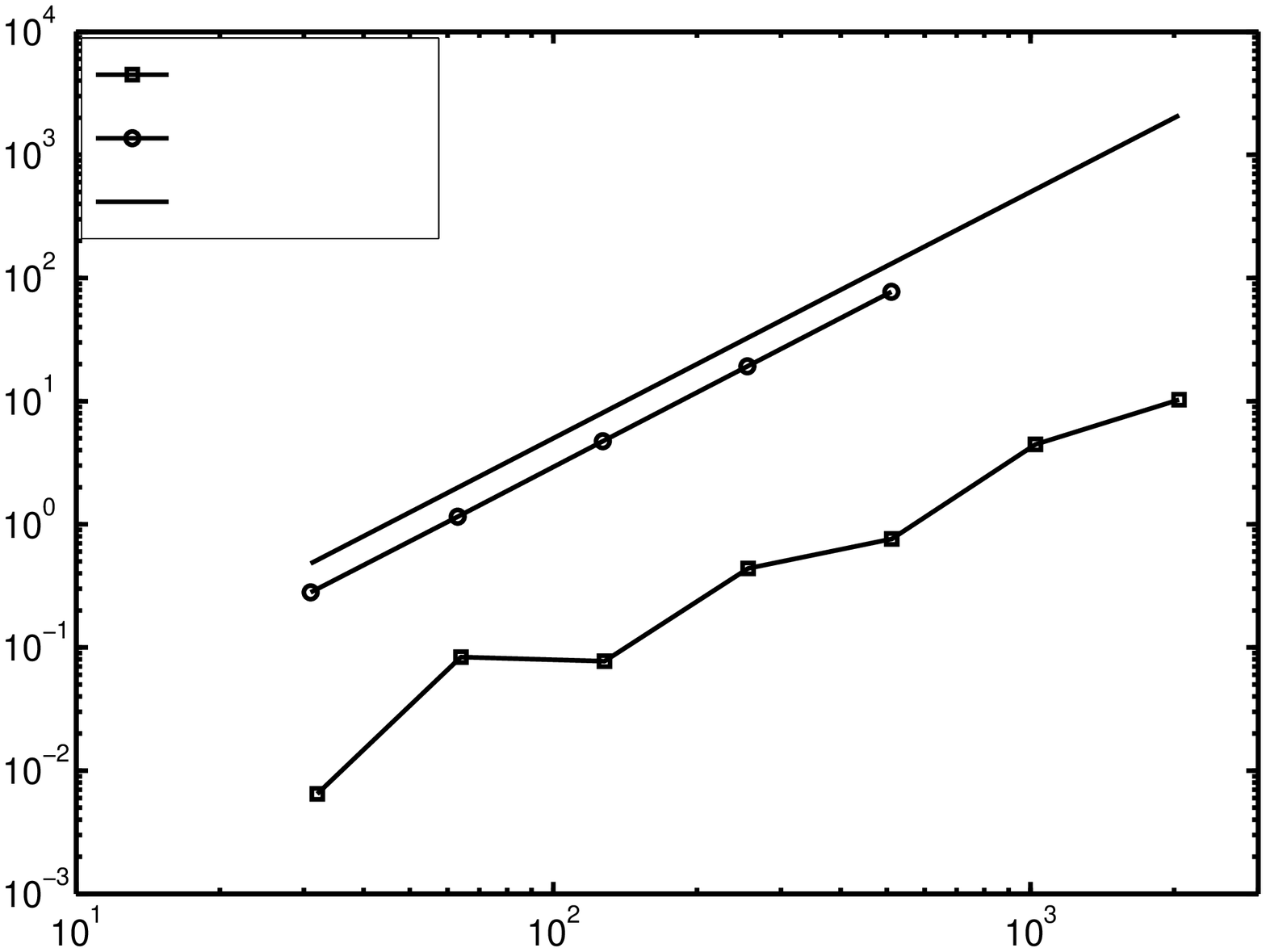}}
\put(75,5){\includegraphics[width=85mm]{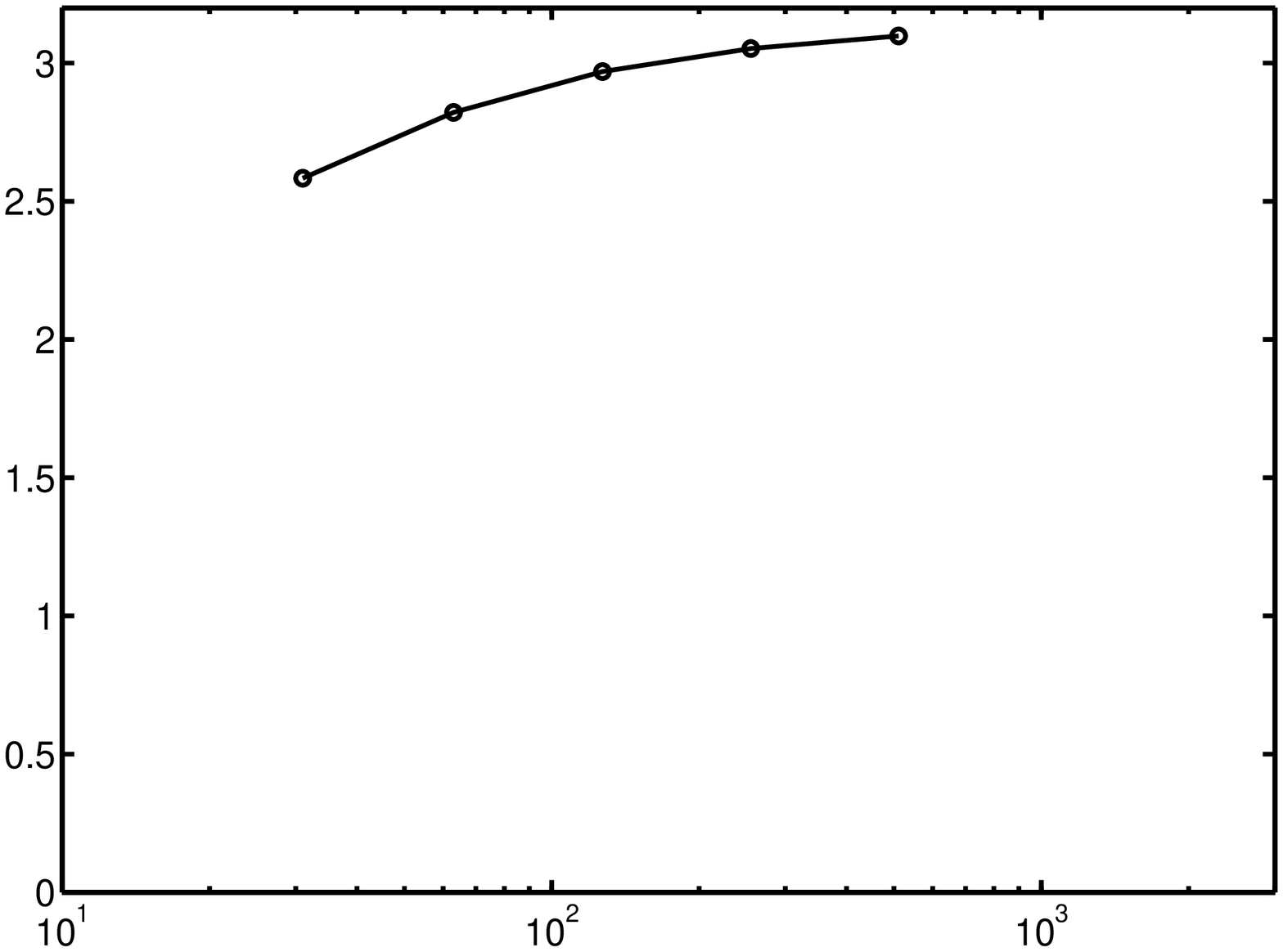}}
\put(38,5){$n$}
\put(118,5){$n$}
\put(37,00){(a)}
\put(117,00){(b)}
\put(-1,25){\rotatebox{90}{\footnotesize  Time in Seconds}}
\put(80,30){\rotatebox{90}{\footnotesize  $M/N_{\rm source}$}}
\put(12,56.0){\scriptsize  FFT}
\put(12,53){\scriptsize   FMM}
\put(12,49){\scriptsize   $10^{-4} \,n^2$}
\end{picture}
\caption{\label{fig:dense} Computational profile for a dense source distribution in a $n \times n$  domain.  Computational times using
the lattice FMM and FFT are reported (a).  The memory $M$ (in KB) per source point ($M/N_{\rm source}$) used in storing the precomputed
information for the lattice FMM are reported (b).}
\end{figure}

In the first experiment, we suppose that every node in the lattice is loaded, see Figure \ref{fig:geometries}(a),
so that $N_{\rm source} = N_{\rm domain} = n^{2}$. In this case, we expect
$$
T_{\rm fft} \sim n^{2}\log(n),
\qquad\mbox{and}\qquad
T_{\rm FMM} \sim n^{2},
$$
and the purpose of the numerical experiment is simply to see how the constants of
proportionality compare.  Figure \ref{fig:dense}(a) provides the answer. We see
that the FMM is slower by roughly two order of magnitude.

In the next three experiments, we suppose that $f$ is only sparsely supported in the domain $\Omega$, so that $N_{\rm source} \ll N_{\rm domain}$. In this case, we expect
$$
T_{\rm fft} \sim n^{2}\log(n),
\qquad\mbox{and}\qquad
T_{\rm FMM} \sim N_{\rm source}.
$$

\begin{figure}[ht]
\centering
\setlength{\unitlength}{1mm}
\begin{picture}(160,85)
\put(-5,5){\includegraphics[width=85mm]{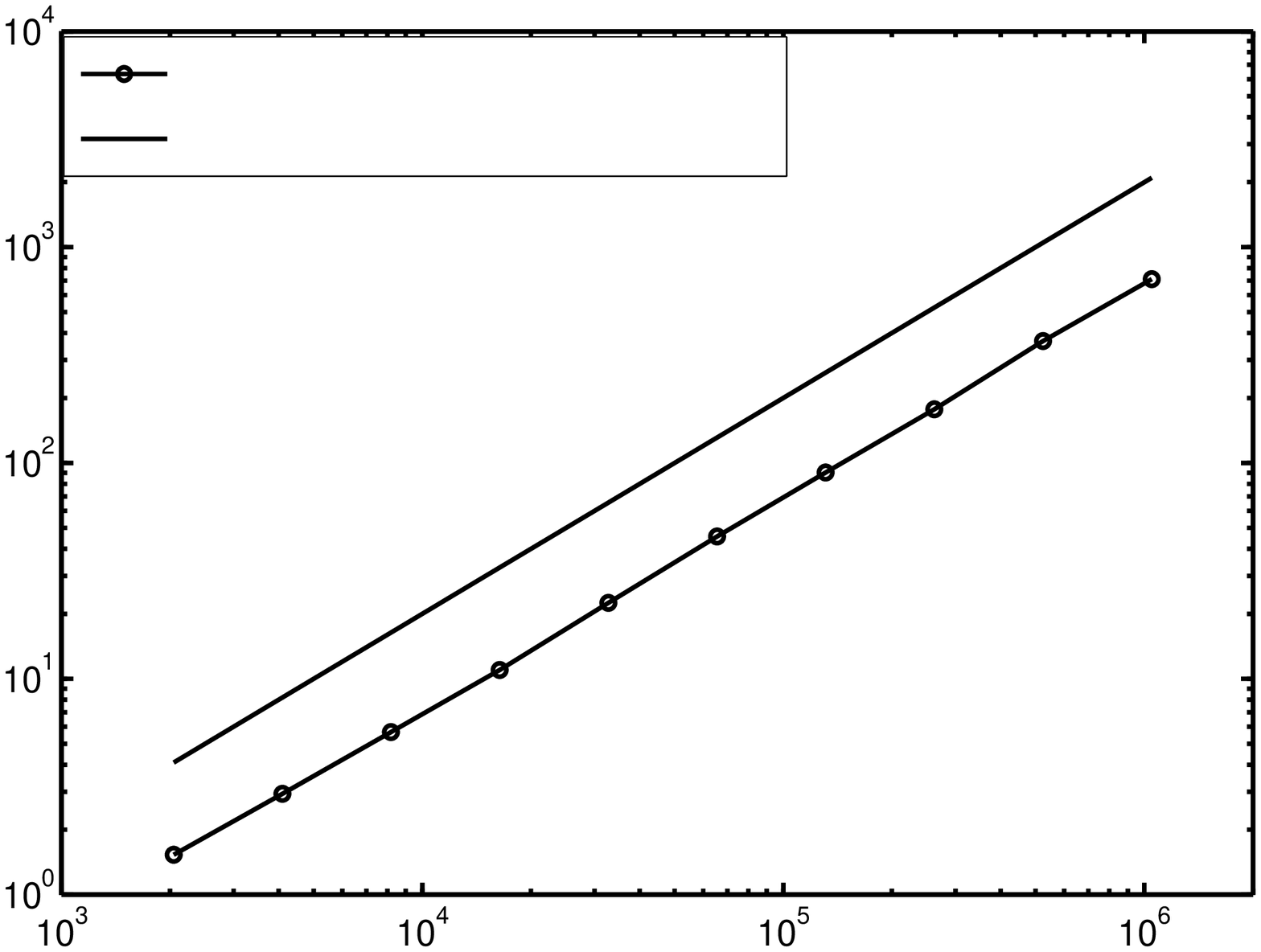}}
\put(75,5){\includegraphics[width=85mm]{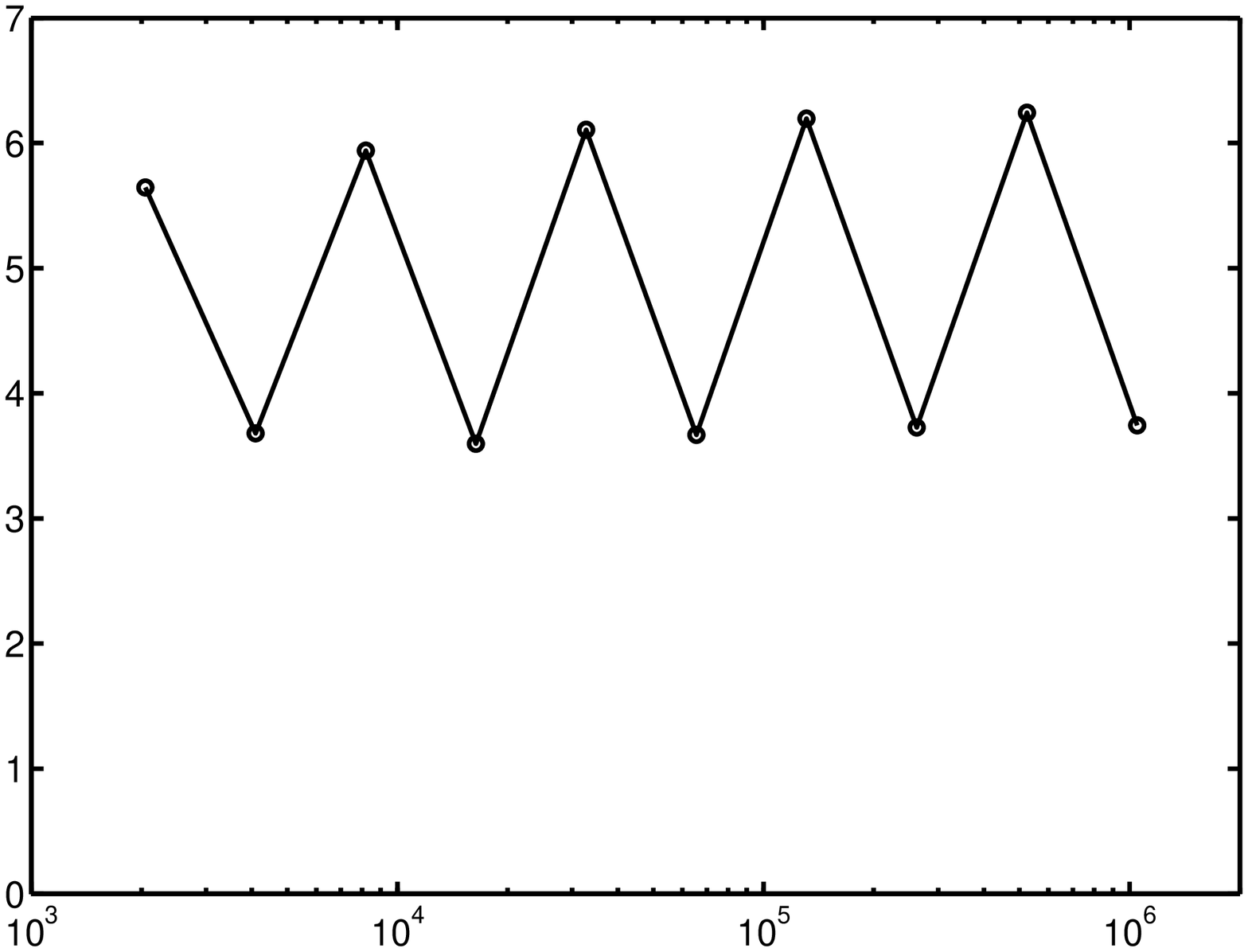}}
\put(38,5){$n$}
\put(118,5){$n$}
\put(37,00){(a)}
\put(117,00){(b)}
\put(-1,25){\rotatebox{90}{\footnotesize  Time in Seconds}}
\put(80,30){\rotatebox{90}{\footnotesize  $M/N_{\rm source}$}}
\put(14,56.5){\scriptsize $n$ uniformly distributed}
\put(14,53){\scriptsize  $10^{-3}\, n$}
\end{picture}
\caption{\label{fig:rand}Computational profile for a $n$ source points distributed via uniform random distribution in a $n \times n$  domain.
Computational times using the lattice FMM are reported (a).  The memory $M$ (in KB) per source point ($M/N_{\rm source}$) used in storing the precomputed
information for the lattice FMM are reported (b).}
\end{figure}

In the second experiment, we suppose that $n$ loads are
distributed according to a uniform random distribution throughout the domain,
see Figure \ref{fig:geometries}(b).
Figure \ref{fig:rand}(a) provides the measured times. It confirms our expectation
that $T_{\rm FMM}$ does not depend on $N_{\rm domain}$, and indeed, that the FMM can handle
a situation with $n = 10^{6}$ loaded nodes in a domain involving $N_{\rm domain} = 10^{12}$
lattice nodes.  Figure \ref{fig:rand}(b) illustrates the memory (in KB) per source point ($M/N_{\rm source}$) required for storing the
pre-computation information.  It confirms our expectation that the memory (in KB) required for storing the pre-computation information depends
linearly with respect to $N_{\rm source}$.

\begin{figure}[ht]
\centering
\setlength{\unitlength}{1mm}
\begin{picture}(160,85)
\put(-5,5){\includegraphics[width=85mm]{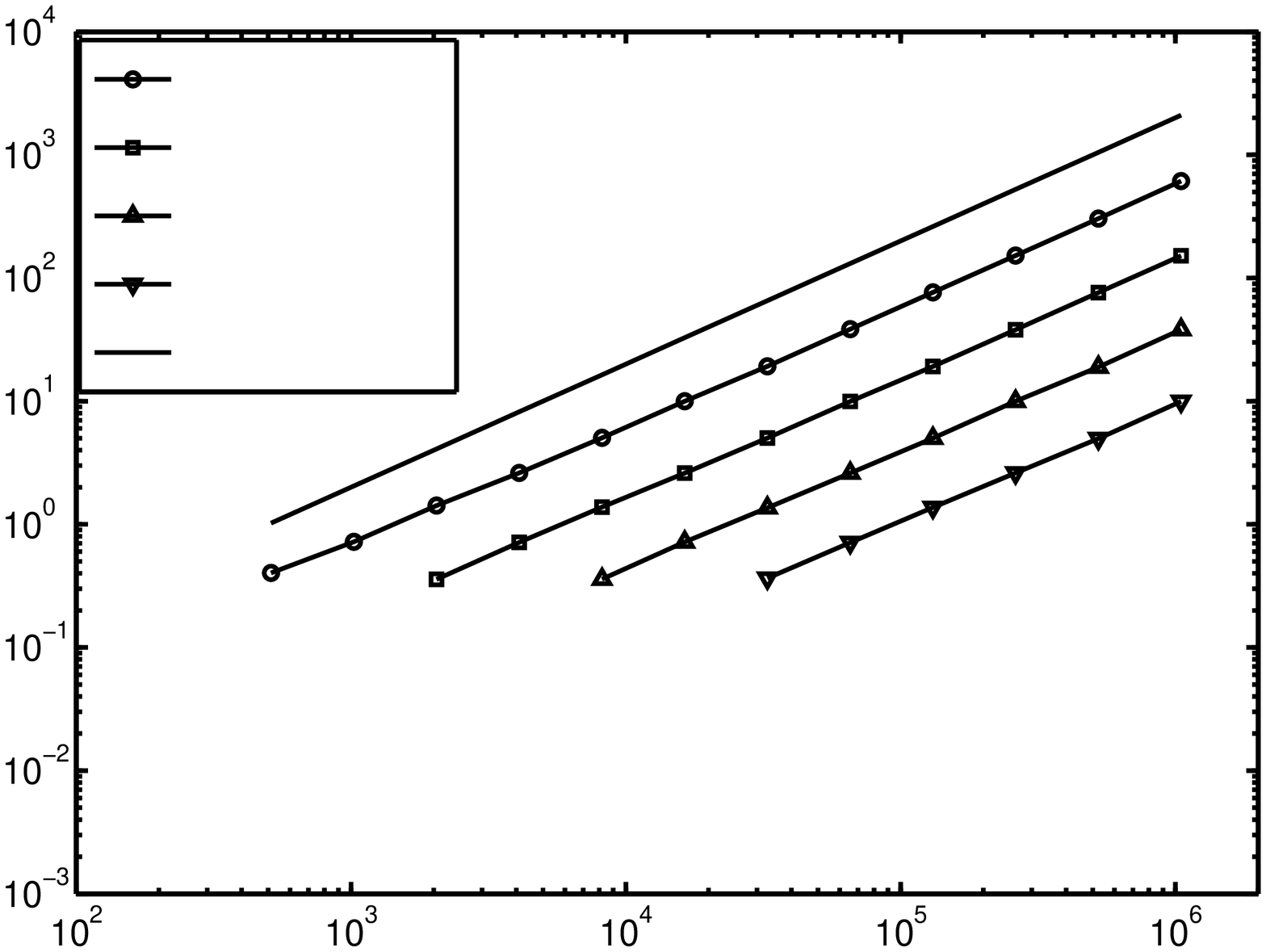}}
\put(75,5){\includegraphics[width=85mm]{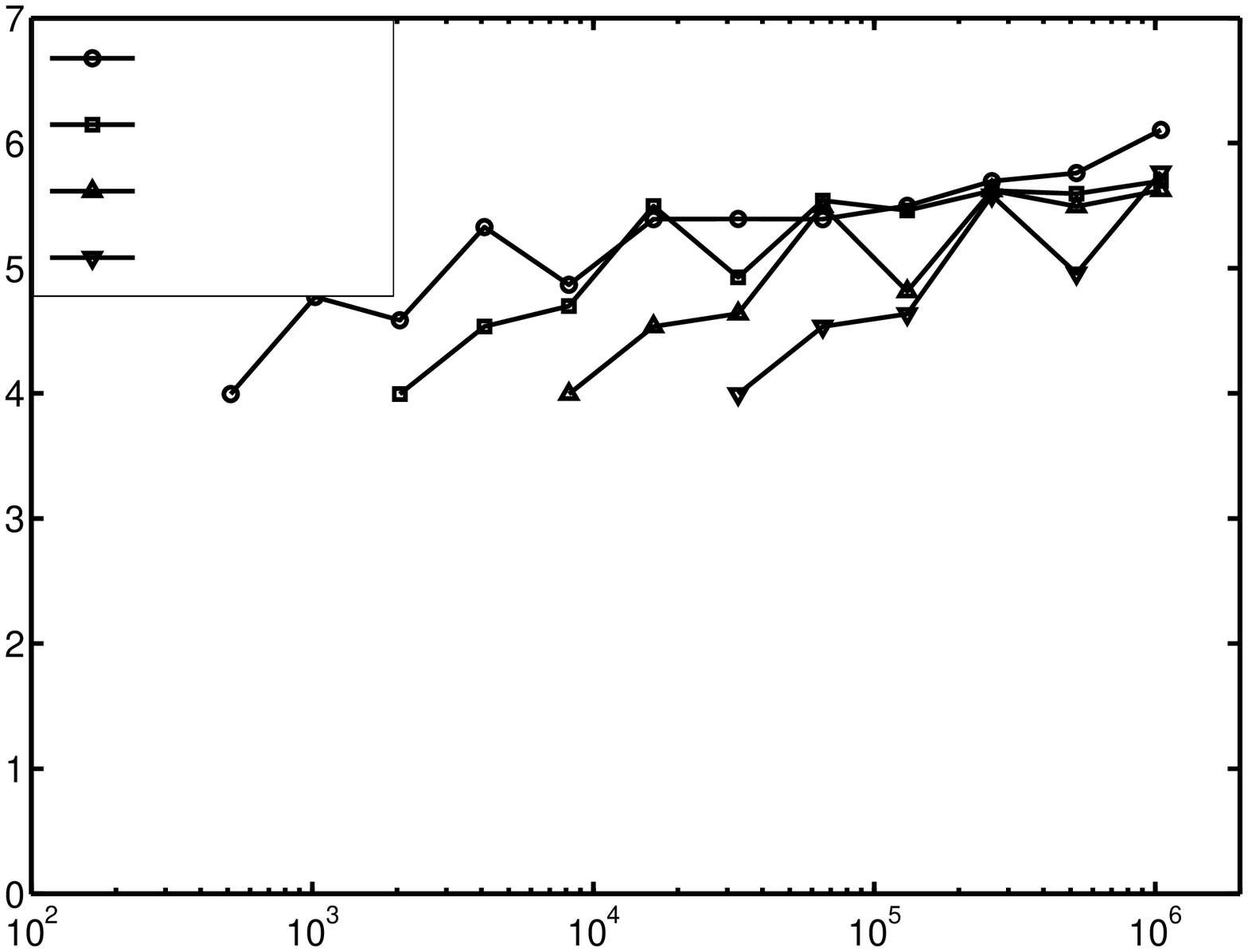}}
\put(38,5){$n$}
\put(118,5){$n$}
\put(37,00){(a)}
\put(117,00){(b)}
\put(-1,25){\rotatebox{90}{\footnotesize  Time in Seconds}}
\put(80,30){\rotatebox{90}{\footnotesize  $M/N_{\rm source}$}}
\put(13,56){\scriptsize   $\alpha = 1$}
\put(13,52.5){\scriptsize   $\alpha = 1/4$}
\put(13,48.5){\scriptsize   $\alpha = 1/16$}
\put(13,45){\scriptsize   $\alpha = 1/64$}
\put(13,41){\scriptsize   $10^{-3}\, n$}
\put(94,56){\scriptsize   $\alpha = 1$}
\put(94,53){\scriptsize   $\alpha = 1/4$}
\put(94,49.5){\scriptsize   $\alpha = 1/16$}
\put(94,46){\scriptsize   $\alpha = 1/64$}
\end{picture}
\caption{\label{fig:circ}Computational profile for $\alpha n$ sources lie on a circle embedded in an $n \times n$  domain.
Computational times using the lattice FMM are reported (a).  The memory $M$ (in KB) per source point ($M/N_{\rm source}$) used in storing the precomputed
information for the lattice FMM are reported (b).}
\end{figure}

In the third experiment, we distribute the load on a circle inscribed in the square $\Omega$, see Figure \ref{fig:geometries}(c), in such a way
that $N_{\rm source} = \alpha\,n$ nodes are loaded, for $\alpha = 1,\ 1/4,\ 1/16, \ 1/64$. Figure \ref{fig:circ}(a)  provides the time measurements
and again confirms our expectation that the $T_{\rm FMM}$ is not dependent on $N_{\rm domain}$.

\begin{figure}[ht]
\centering
\setlength{\unitlength}{1mm}
\begin{picture}(160,90)
\put(-5,5){\includegraphics[width=85mm]{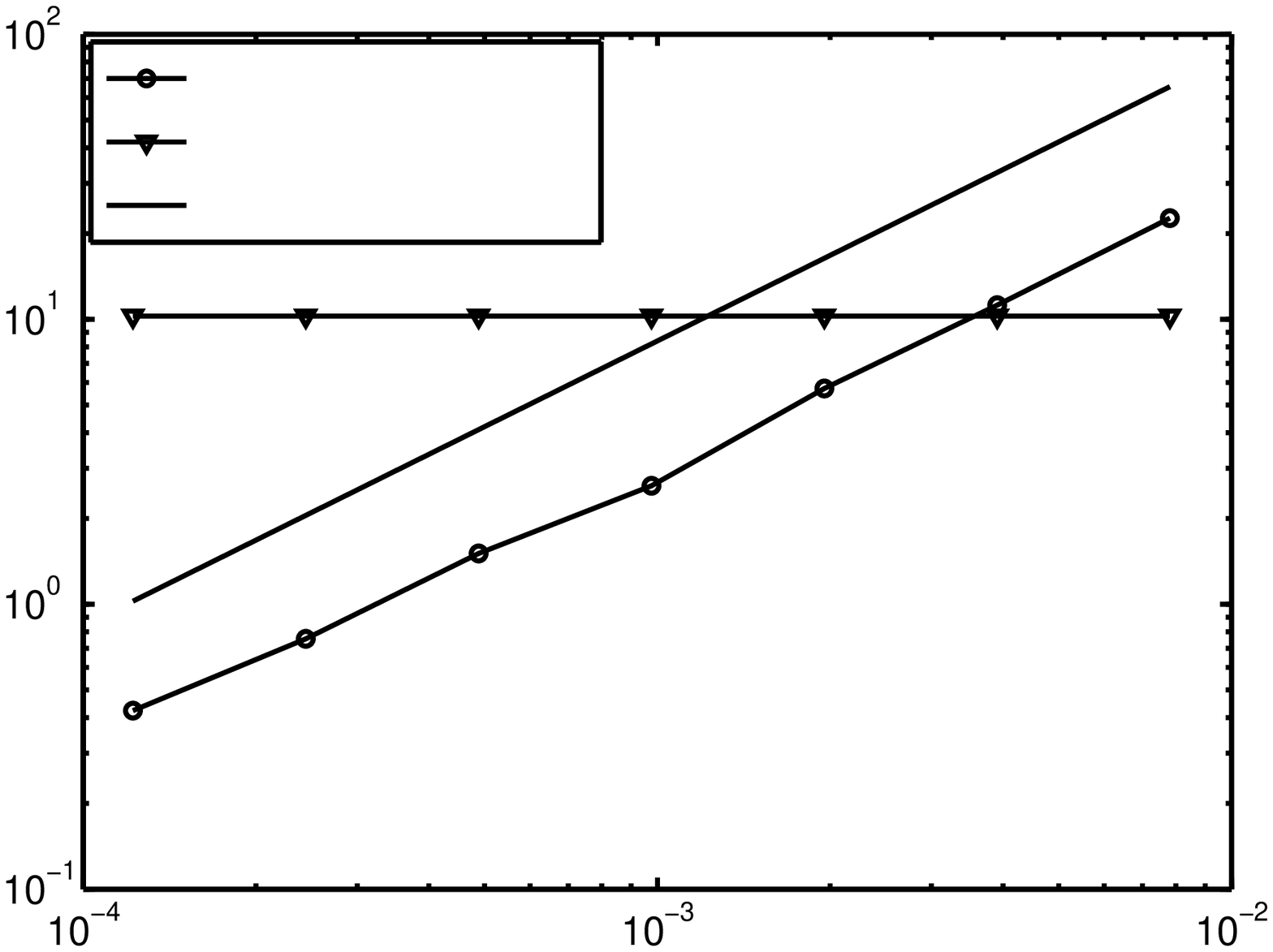}}
\put(75,5){\includegraphics[width=85mm]{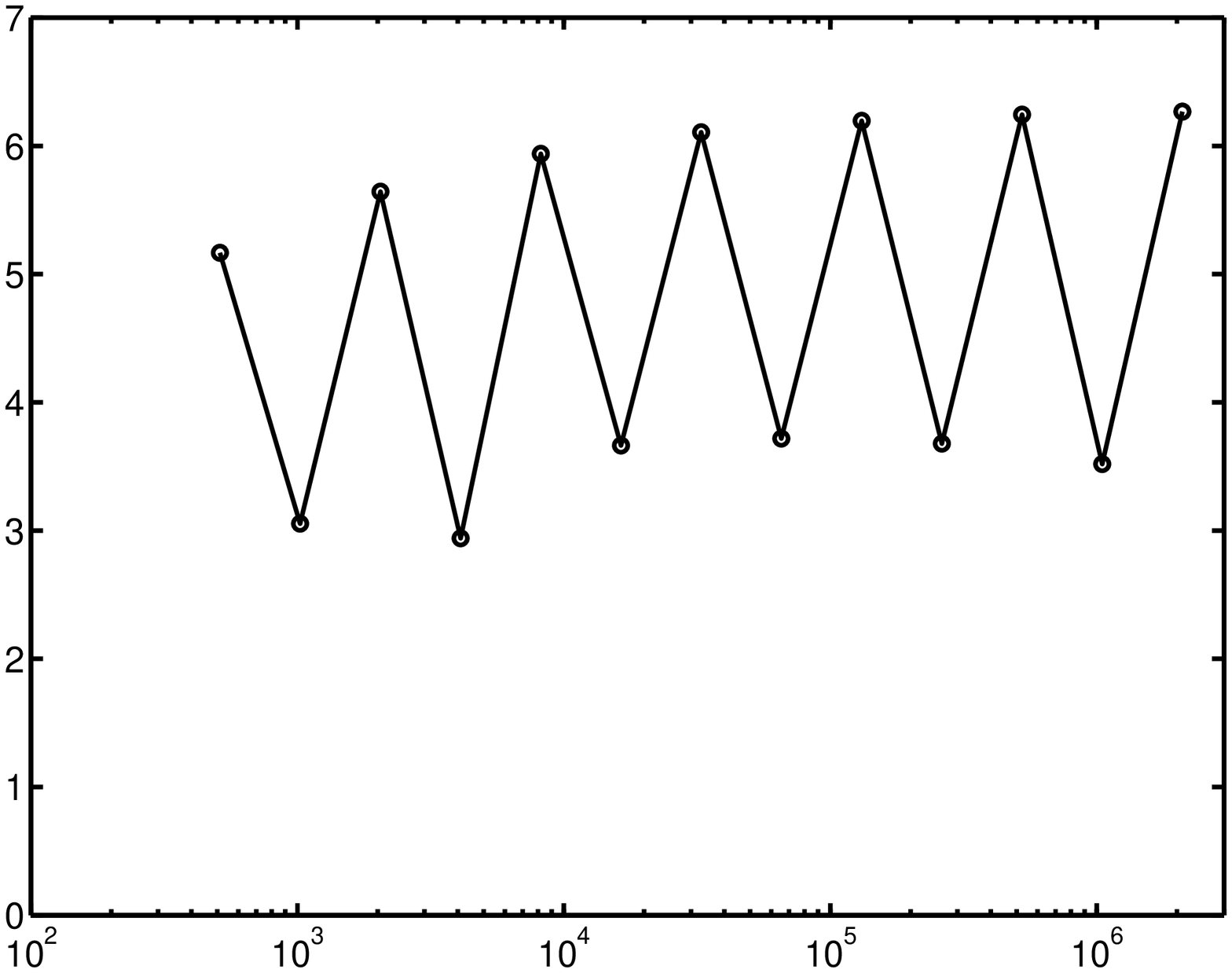}}
\put(28,5){$N_{\rm source}/N_{\rm domain}$}
\put(115,5){$N_{\rm source}$}
\put(37,00){(a)}
\put(117,00){(b)}
\put(-1,25){\rotatebox{90}{\footnotesize  Time in Seconds}}
\put(80,30){\rotatebox{90}{\footnotesize  $M/N_{\rm source}$}}
\put(13,57){\scriptsize   Random sources}
\put(13,53){\scriptsize   FFT}
\put(13,50){\scriptsize   $10^{-3} \,N_{\rm source}$}
\end{picture}
\caption{\label{fig:fixed}Computational profile for a fixed  $2048 \times 2048$ lattice domain.
Computational times using the lattice FMM and the FFT are reported (a).  The memory $M$ (in KB) per source point ($M/N_{\rm source}$) used in storing the precomputed
information for the lattice FMM are reported (b).}
\end{figure}

In the final experiment, we fix the domain to be sized $2048 \times 2048$ and increase the number of body loads distributed according to a
uniform distribution.
Figure \ref{fig:fixed}(a) provides the time measurements in comparison with the FFT.  It illustrates that for sources occupying less than $0.39 \%$
of the domain (corresponding to $16,384$ sources) the lattice FMM is the faster method.

\clearpage

\section{Extensions}

\subsection{Lattices with inclusions}
\label{sec:extend_inclusions}

\begin{figure}
\centering
\includegraphics[height=4.5cm]{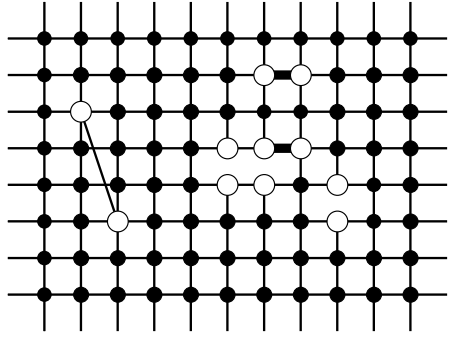}
\caption{A piece of an infinite lattice with some deviations from perfect periodicity.
One bar has been added, three bars have been removed, and two bars have been strengthened.
The set $\Omega_{\rm inc}$ of effected nodes has 11 elements, which are marked with
white circles.}
\label{fig:inclusions}
\end{figure}

The lattice FMM described in this paper can be used to handle many lattices
featuring local deviations from perfect periodicity due to, \textit{e.g.},
missing bars or lattice inclusions, see Figure \ref{fig:inclusions}. To
illustrate, let us consider the equation
\begin{equation}
\label{eq:pert}
[(\mtx{A} + \mtx{B})u](\vct{m}) = 0,\qquad \vct{m} \in \mathbb{Z}^{2},
\end{equation}
where $\mtx{A}$ is the discrete Laplace operator (\ref{eq:def_dL}),
and where $\mtx{B}$ is an operator corresponding to some local modifications
to the lattice, such as those illustrated in Figure \ref{fig:inclusions}.
A typical far-field condition of interest
is to require that the potential $u$ approaches a linear function, say
$v(\vct{m}) = c_{1}m_{1} + c_{2}m_{2}$ where $c_{1}$ and $c_{2}$ are given constants, at infinity. Formally, we require
\begin{equation}
\label{eq:pert_decay}
\lim_{|\vct{m}|\rightarrow\infty} \bigl| u(\vct{m}) - v(\vct{m})\bigr| = 0.
\end{equation}
(In what follows, any function $v$ satisfying $\mtx{A}v = 0$ would work.)
When we have access to the solution operator $\mtx{S}$, defined via
$$
[\mtx{S}u](\vct{m}) = \sum_{\vct{n} \in \mathbb{Z}^{2}}\phi(\vct{m}-\vct{n})\,u(\vct{n}),
$$
then the equation (\ref{eq:pert}), which is a sparse equation defined on the infinite
lattice $\mathbb{Z}^{2}$, can be transformed to a dense equation defined only on the
set $\Omega_{\rm inc}$ of nodes which connect to bars whose conductivity has been changed,
or to which new bars have been connected.

To execute the reduction of (\ref{eq:pert}) to an equation defined on $\Omega_{\rm inc}$, set
$$
u = v + w,
$$
and observe that since $\mtx{A}v = 0$ and $(\mtx{A}+\mtx{B})u=0$, the new unknown $w$ must satisfy
\begin{equation}
\label{eq:eq_for_w}
(\mtx{A} + \mtx{B})\,w = -\mtx{B}\,v.
\end{equation}
Since $\lim_{|\vct{m}|\rightarrow\infty} \bigl|w(\vct{m})\bigr| = 0$, we have $\mtx{S}\,\mtx{A}\,w = w$,
and so application of $\mtx{B}\mtx{S}$ to (\ref{eq:eq_for_w}) results in the equation
\begin{equation}
\label{eq:pert_fred}
\mu + \mtx{B}\,\mtx{S}\,\mu = -\mtx{B}\,\mtx{S}\,\mtx{B}\,v,
\end{equation}
where
$$
\mu = \mtx{B}\,w.
$$
The key observation is now that $\mtx{B}$ is a local operator on $\Omega_{\rm inc}$, so
equation (\ref{eq:pert_fred}) can be restricted to $\Omega_{\rm inc}$ to obtain a closed equation for $\mu$.
Once $\mu$ has been solved from this equation, the original potential $u$ is recovered via
$$
u = v - \mtx{S}\,(\mtx{B}\,v - \mu).
$$

It has been demonstrated \cite{gillman} that iterative solvers in many situations
converge rapidly when used to solve a system such as (\ref{eq:pert_fred}). The
dense coefficient matrix can rapidly be
applied to a vector since $\mtx{B}$ is a sparse operator, and since $\mtx{S}$
is amenable to the lattice FMM described in this paper.

\subsection{Finite lattice problems involving boundary conditions}
\label{sec:extend_finite_lattices}
The lattice FMM described in this paper can also be used in a solution technique for boundary value problems.
To illustrate, let us consider the Dirichlet problem

\begin{equation}
\label{eq:Dirichlet}
\left\{\begin{split}
[\mtx{A}\,u](\vct{m}) &= 0,\hspace{14mm} \vct{m} \in \Omega,\\
u(\vct{m}) &= g(\vct{m}),\qquad \vct{m} \in \Gamma.
\end{split}\right.
\end{equation}

When we have access to the solution operator $\mtx{D}$, defined

\begin{equation}\label{eq:Deq} [\mtx{D}u](\vct{m}) = \partial_{\vct{n}}\phi(\vct{m}-\vct{n}) q(n)\end{equation}
where the subscript $\vct{n}$ in $\partial_{\vct{n}}$ simply indicates that the difference operator is
acting on the variable $\vct{n}$, then the equation (\ref{eq:Dirichlet}) can be reduced to a problem defined on $\Gamma$.

An equation for $q$ is obtained by simply restricting
(\ref{eq:Deq}) to $\Gamma$:
\begin{equation}
\label{eq:wellcond}
\sum_{n \in \Gamma} d(\vct{m},\vct{n})\,q(\vct{n}) = g(\vct{m}),
\qquad \vct{m} \in \Gamma.
\end{equation}

It has been shown in \cite{gillman} that there exist an $O(N_\Gamma)$ inversion scheme for solving (\ref{eq:wellcond}).  The operator $\mtx{D}$ can
 be applied rapidly to $q$ using the lattice FMM described in this paper.

Also in \cite{gillman}, the technique described in this section and in section \ref{sec:extend_inclusions} has been put together to create a
fast solver for finite lattice problems involving boundary conditions and inclusions.  A key tool for the fast solver is the lattice FMM described
in this paper.

\section{Concluding remarks}
The paper presents a kernel independent FMM for solving Poisson problems defined on the integer lattice $\mathbb{Z}^2$.  For simplicity of
presentation, we focused on equations involving the discrete Laplace operator.  Techniques for evaluating the corresponding lattice fundamental
solutions are presented.  The complexity of the proposed method is $O(N_{\rm source})$ where $N_{\rm source}$ is the number of locations in
$\mathbb{Z}^2$ subjected to body loads.

Numerical experiments demonstrate that for problems where the body loads are sparsely distributed
in a computational box the proposed method is faster and more robust than the FFT.  For instance, it was
 demonstrated that using a standard desktop PC, a lattice Poisson equation on a lattice with $N_{\rm domain} = 10^{12}$
nodes, of which $N_{\rm source} = 10^6$ were loaded, was solved to ten digits of accuracy in
three minutes. It should be noted that this problem is about six orders of magnitude
larger than the largest Poisson problem that can be handled via the FFT.  Also, it was demonstrated for a lattice Poisson problem
in a domain with $N_{\rm domain} = 4,194,304$ nodes, the lattice FMM is faster than the FFT when the number of loaded points is
less than $N_{\rm source} = 16,384$.

\bibliography{main_bib}
 \bibliographystyle{amsplain}
\end{document}